\documentclass{article}

\usepackage{geometry}
\usepackage{amsthm}
\usepackage{mathptmx}
\usepackage{amssymb}
\usepackage{amsfonts}
\usepackage{tikz}
\usepackage{epstopdf}
\usepackage{authblk}
\usepackage{amsmath,amscd}
\usepackage{float}
\usepackage[caption = false]{subfig}
\usepackage{pgfplots}
\usepackage{tikz}
\usepackage{amscd}
\usepackage{graphicx}
\usepackage{color}
\usepackage{bold-extra}
\usepackage{algorithm2e}
\usepackage{dcolumn}
\usepackage{ulem}  
\usetikzlibrary{arrows.meta}

\DeclareMathOperator*{\argmax}{arg\,max}
\DeclareMathOperator*{\argmin}{arg\,min} 

\newcolumntype{d}[1]{D{.}{.}{#1}}
\newcolumntype{m}[1]{D{\pm}{\pm}{#1}}

\newcommand{\R}{\mathbb{R}}
\newcommand{\bfn}{\boldsymbol{n}}

\begin{document}

\title{Monte Carlo on manifolds: sampling densities and integrating functions}

\author[1]{Emilio Zappa}  
\author[1]{Miranda Holmes-Cerfon}
\author[1]{Jonathan Goodman}

\affil[1]{Courant Institute of Mathematical Sciences, New York University, NY.}

\maketitle

\begin{abstract}
We describe and analyze some Monte Carlo methods for manifolds in Euclidean space defined by 
equality and inequality constraints.
First, we give an MCMC sampler for probability distributions defined by un-normalized
densities on such manifolds.
The sampler uses a specific orthogonal projection to the surface that requires only information about the tangent space to the manifold, obtainable from first derivatives of the constraint functions, hence avoiding the need for curvature information or second derivatives.
Second, we use the sampler to develop a multi-stage algorithm to compute integrals over such manifolds.
We provide single-run error estimates that avoid the need for multiple independent runs.
Computational experiments on various test problems show that the algorithms and error 
estimates work in practice.
The method is applied to compute the entropies of different sticky hard sphere systems. These predict the temperature or interaction energy at which loops of hard sticky spheres become preferable to chains.
\end{abstract}

\section{Introduction} \label{sec:intro}

Many physical and engineering systems involve motion with equality constraints.
For example, in physics and chemistry, bonds between atoms or colloidal particles may be modeled as having a fixed length \cite{miranda}, while in robotics, joints or hinges connecting moving components may be written as constraints on distances or angles \cite{Wampler:2011gq}. Constraints also arise in statistics on the set of parameters in parametric models \cite{diaconis}. In such problems the space of accessible configurations is lower-dimensional than the space of variables which describe the system, often forming a manifold embedded in the full configuration space. 
One may then be interested in sampling a probability distribution defined on the manifold, 
or calculating an integral over the manifold such as its volume. 

This paper presents Monte Carlo methods for computing a $d$-dimensional integral over a 
$d$-dimensional connected manifold $M$ embedded in a $d_a$-dimensional Euclidean space, 
with $d_a\geq d$.
The manifold is defined by sets of constraints.
We first describe an MCMC sampler for manifolds defined by constraints in this way.
Then we describe a multi-phase procedure that uses the sampler to estimate integrals.
It should go without saying (but rarely does) that deterministic algorithms to compute integrals are
impractical except in very low dimensions or for special problems  (e.g.\cite{simonovits}).

Our  MCMC sampler has much in common with other algorithms for constraint manifolds, including
other samplers, optimizers, differential algebraic equation (DAE) solvers 
(see, e.g., \cite{hairer}), etc.  
A move from a point $x \in M$ starts with a move in the tangent space at $x$, which is $v \in T_x$.
This is followed by a projection back to $y \in M$, which we write as $y = x + v + w$.
The sampler simplifies if 
we require $w \perp T_x$, so the projection is 
perpendicular to the tangent space \emph{at the original point}. 
This idea was first suggested in \cite{lelievre2012}\footnote{We learned of this work after we submitted the original version of this paper for publication.}.
This is different from other natural choices, such as choosing $y \in M$
to minimize the projection distance $\left\| y-(x+v)\right\|$. 
The choice of $w\perp T_x$ makes the sampler easy to implement as it requires only first derivatives of the constraint surfaces.
Neither solving for $y$ using Newton's method nor the Metropolis Hastings detailed balance condition require any higher derivatives.
Other surface sampling methods that we know of require second derivatives, 
such as the explicit parameterization method of \cite{diaconis}, the Hamiltonian method
such as of \cite{brubaker}, 
geodesic methods such as \cite{girolami,leimkuhler2016}, 
and discretized SDE methods such as \cite{ciccotti}.
It is common in our computational experiments that there is no $w \perp T_x$ so that $x+v+w \in M$, or that one or more $w$ values exist 
but the Newton solver fails to find one.
Section \ref{sec:MCMC} explains how we preserve detailed balance even in these situations.

Our goal is to compute integrals of the form 
\begin{equation}  \label{eq:Z}
      Z =   \int_M f(x) \,\sigma(dx) \; ,
\end{equation}
where $\sigma(dx)$ is the $d-$dimensional surface area measure (Hausdorff measure) and 
$f$ is a positive smooth function.
Writing $B \subset \R^{d_a}$ for the ball centered at a point $x_0\in M$, we also estimate
certain integrals of the form 
\begin{equation}  \label{eq:ZB}
      Z_B =   \int_{M\cap B} f(x) \,\sigma(dx) \; .
\end{equation}
For each $B$, we define a probability distribution on $M\cap B$:
\begin{equation}  \label{eq:rhoB}
        \rho_B(dx) = \frac{1}{Z_B} f(x) \sigma(dx) \; .
\end{equation}
The sampling algorithm outlined above, and described in more detail in Section \ref{sec:MCMC}, 
allows us to draw samples from $\rho_B$.

To calculate integrals as in \eqref{eq:Z} we use a multi-phase strategy related to nested sampling and thermodynamic integration (e.g.\cite{frenkel84}).
Similar strategies were applied to volume estimation problems by the Lovasz school, 
see, e.g., \cite{vempala}.
Let $B_0 \supset B_1 \supset \cdots \supset B_k$ be a nested contracting family of balls with 
a common center.
With an abuse of notation, we write $Z_i$ for $Z_{B_i}$.
We use MCMC sampling in $M\cap B_{i}$ to estimate the ratio
\begin{equation}   \label{eq:Ri}
       R_i = \frac{Z_{i}}{Z_{i+1}} \; .
\end{equation}
We choose the last ball $B_k$ small enough that we may estimate $Z_k$ by direct Monte-Carlo integration.
We choose $B_0$ large enough that $M \subset B_0$ 
so $Z = Z_0$.
Our estimate of $Z$ is 
\begin{equation}  \label{eq:Zhat}
       \widehat{Z} = \widehat{Z}_k \, \prod_{i=0}^{k-1} \widehat{R}_i \; ,
\end{equation}
where hatted quantities such as $\widehat{Z}$ are Monte Carlo estimates.
Section \ref{sec:vol} describes the procedure in more detail.

Error estimates should be a part of any Monte Carlo computation.
Section \ref{sec:errors} describes our procedure for estimating 
\[
        \sigma^2_{\widehat{Z}} = \mbox{var} \!\left( \widehat{Z} \right) \; ,
\]
a procedure adapted from  \cite{hou}.
If we can estimate $\sigma^2_{\widehat{R}_i}$, and if these are small enough, then we can combine
them to estimate $\sigma^2_{\widehat{Z}}$.
We explore how $\sigma_{\widehat{Z}}^2$ depends on parameters of the algorithm via several simplified model problems, in section \ref{sec:toy}. 
We specifically analyze how the error depends on the amount by which we shrink the balls for each ratio, and make predictions that could lead to good systematic choices of this parameter.

We illustrate our methods on a number of examples.
We apply the sampler to the surface of a torus and cone in 3D, as well as the special orthogonal group $SO(n)$ for $n=11$, which is a manifold of dimension $55$ in an ambient space of dimension 121 (section \ref{sec:examples}) . All of these examples have certain marginal distributions which are known analytically, so we verify the sampler is working correctly. 
We verify the integration algorithm is correct by calculating the surface areas of a torus and $SO(n)$ for $2\leq n \leq 7$ (sections \ref{sec:torusvol}, \ref{sec:SO(n)vol}.) 

Finally, we apply our algorithm to study clusters of sticky spheres (section \ref{sec:spheres}), a model that is studied partly to understand how colloidal particles self-assemble into various configurations \cite{colloids2,miranda,mhc2017}.
A simple question in that direction concerns a linear chain of $n$ spheres -- whether it is more likely to be closed in a loop or open in a chain.
According to equilibrium statistical mechanics, the ratio of the probabilities of being open versus closed depends on a parameter ( the ``sticky parameter'') characterizing the binding energy and temperature, and the entropies, which  are surface areas of the kind we are computing. 
We compute the entropies for open and closed chains of lengths 4--10. This data can be used to predict the relative probabilities of being in a loop or a chain, if the sticky parameter is known, but if it is not known our calculations can be used to measure it by comparing with data.

\paragraph{Notation and assumptions.} Throughout the paper, we will denote by $M$ a $d$-dimensional 
connected manifold embedded in an ambient space $\R^{d_a}$, which is defined by equality and 
inequality constraints.
There are $m$ equality constraints, $q_i(x) = 0$, $i = 1, \ldots, m$,
where the $q_i$ are smooth functions $q_i : \R^{d_a} \rightarrow \R$.
There are $l$ inequality constraints of the form $h_j(x) > 0$, for $j =1, \ldots, l$.
The manifold $M$ is 
\begin{equation}\label{eq:manifold}
M = \left\{ x \in \R^{d_a} : q_i(x) = 0, \; i=1, \ldots, m, \; h_j(x) > 0, \; j = 1, \ldots, l \right\}.
\end{equation}
If $M$ has several connected components, then the sampler may sample only one
component, or it may hop between several components.
It is hard to distinguish these possibilities computationally. 

We let $Q_x$ be the matrix whose columns are the gradients $\{ \nabla q_i (x) \}_{i=1}^m$.
This makes $Q_x$ the transpose of the Jacobian of the overall constraint function 
$q: \R^{d_a} \to R^m$.
The entries of $Q_x$ are 
\begin{equation}\label{eq:Qx}
\left( Q_x \right)_{ij}  = \frac{\partial q_i(x)}{\partial x_j} \; .
\end{equation}
We assume that $Q_x$ has full rank $m$ everywhere on $M$. 
By the implicit function theorem, this implies that the dimension of $M$ is $d = d_a-m$ and that the tangent space $T_x\equiv T_xM$ at a point $x \in M$ is well defined. In this case $\{ \nabla q_i (x) \}_{i=1}^m$ form a basis of the orthogonal space $T^\perp_x \equiv T_xM^{\perp}$. 
Note that $M$ inherits the metric from the ambient space $\R^{d_a}$ by restriction. 
The corresponding $d-$dimensional volume element is $d-$dimensional Hausdorff measure, which we
denote by $\sigma(dx)$.

\section{Sampling on a constraint manifold}  \label{sec:MCMC}

\begin{figure}
\noindent\fbox{%
    \parbox{\textwidth}{%
\begin{center}
    MCMC surface sampling algorithm
\end{center}
\fontsize{8}{11}\selectfont
\begin{tabbing}
\hspace{2em}\=\hspace{2em}\=\hspace{2em}\=\hspace{2em}\=\hspace{2em}\=\hspace{2em}\= \hspace{16em}\= \\
{\bf Begin:} $x = X_n$, $f(x)$ is known $\left\{\right.$ \\
    \>{\bf Proposal:} generate $y \in M$ or declare failure\\
    \>    \> Calculate $Q_x$ using (\ref{eq:Qx})\\
    \>    \> Find orthonormal bases for $T_x$ and $T_x^{\perp}$ using the QR decomposition of $Q_x$\\
    \>    \> Generate $v \in T_x$ with $v \sim p(v)$ using the orthonormal basis of $T_x$\\
    \>    \> {\bf Projection} to $M$:  \\
        \>    \>   \>    $\left[a,\texttt{flag}\right]$ = {\bf project}$(q,x+v,Q_x)$ \>\>\>\>
// solve $q(x+v+Q_xa) = 0$ for $a$ \\
 \>\>  \> {\bf if}$\;\; \texttt{flag} == 0 \;\; \left\{ \right.$\\
  \>\>    \>     \> $X_{n+1} = x$;  {\bf done} 
                   \>\>\>  // projection failure: reject,  back to Begin.\\
   \>\>    \>  $\left. \;\;\;\;\right\}$ \\
  \>    \>   $y = x + v + Q_x a$
                  \>\>\>\> // successful projection \\
      \>     {\bf Inequality} check: \\
      \>\>{\bf if} $h_i(y) < 0$ for some $y$ $\left\{\right.$ \\
      \>    \>\>   $X_{n+1} = x$;  {\bf done}  \>  \>\>\>// inequality failure: reject; back to Begin.\\
      \>   \>  $\left. \;\;\;\;\right\}$  \\
    \>{\bf Metropolis-Hastings acceptance/rejection step:} \\
    \>    \> Generate $Q_y$ and find orthonormal bases for $T_y$ and $T_y^{\perp}$\\
    \>    \> Find $v^{\prime} \in T_y$ and $w^{\prime} \in T_y^{\perp}$ so that $x = y + v^{\prime} + w^{\prime}$\\
    \>    \> $\displaystyle P_a = \min\!\left( 1, \frac{f(y)p(v^{\prime})}{f(x) p(v)}\right)$
                     \>\>\>\>\> // acceptance probability \\
    \>    \> Generate $U \sim \mbox{unif}(0,1)$ \\    
    \>    \> {\bf if} $\;\; U > P_a \;\; \left\{ \right.$   
                   \>\>\>\>\> // $\mbox{Pr}(\mbox{rej}) = 1-P_a$ \\ 
    \>    \>    \> $X_{n+1} = x$;  {\bf done} 
                  \>\>\>\>  // rejection; back to Begin.\\
    \>    \>    $\left. \;\;\;\;\right\}$ \\ 
    \>     {\bf Reverse Projection}, try to recover $x = y + v^{\prime} + Q_ya$\\
     \>    \>      $\left[a,\texttt{flag}\right]$ = {\bf project}$(q,y + v^{\prime},Q_y)$ \>\>\>\>\\
 \>\>  {\bf if}$\;\; \texttt{flag} == 0 \;\; \left\{ \right.$\\
  \>\>  \>$X_{n+1} = x$;  {\bf done} 
                   \>\>\>\>  // reverse projection failure: reject,  back to Begin.\\
   \>\>    \>  $\left. \;\;\;\;\right\}$ \\

    \>       {\bf Accept move: }  $X_{n+1} = y$
                  \>\>\>\> \> \>// successful recovery; accept the proposal \\    
 $\left.\right\}$\\
 \\
 
$\left[a,\texttt{flag}\right]$ = {\bf project}$(q,z,Q)$
 \>\>\>\>\>\>\>  // Newton's method to solve $q(z+Qa)=0$ for $a$ \\
      \>Initialization: $a = 0$, $i = 0$, $\texttt{flag}=1$\\
     \>{\bf While} $\;\;\left| q(z+Qa)\right| > \mbox{\texttt{\textbf{tol}}}\;\; \left\{ \right.$\\
     \>   \>   Solve $\left(\nabla q(z+Qa)^t Q \right) \Delta a = -q(z+Qa)$ for $\Delta a$
     \>\> \> \>\> // see (\ref{eq:Ja})\\
      \>    \> $a = a + \Delta a$
                 \\
     \>    \> $i = i+1$ \\
   \>    \> {\bf if}$\;\; i >  \mbox{\texttt{\textbf{nmax}}} \;\; \left\{ \right.$\\
    \>    \>     \> $\texttt{flag}=0$;  {\bf return} 
                   \>\>\>\>  // projection failure: reject;  return.\\
   \>    \>   $\left. \;\;\;\;\right\}$ \\
   \>   $\left. \;\;\;\;\right\}$ 
 \end{tabbing}
    }%
}
\caption{Summary of the MCMC algorithm, in pseudocode. 
} \label{fig:MCMC}
\end{figure}

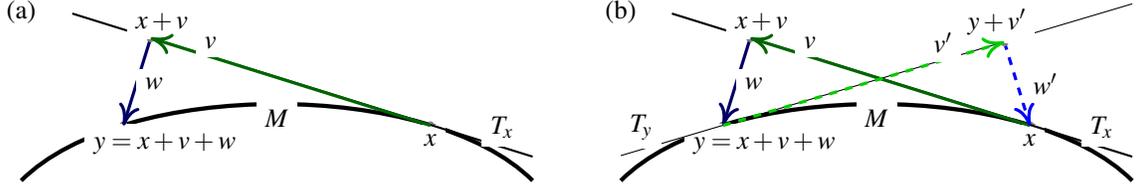
\begin{figure} 
\begin{center}
\begin{tikzpicture}[scale=0.68]
  \draw [ultra thick] (0,0) .. controls (2,2)  and (8,2) .. (10,0);  
  \draw (5,1.23) node[fill=white] {$M$};                             
  
  \filldraw [gray] (8,1.08) circle (2pt);                            
  \draw (8,.8) node[fill=white] {$x$};                               
  \draw [thick] (10,.47) -- (1,3.27);                                
  \draw (9.4,1.05) node[fill=white] {$T_x$};                         
  
  \filldraw [gray] (2,1.08) circle (2pt);                            
  \draw (2.8,.73) node[fill=white] {$y=x+v+w$};                      
  
  \draw [-{>[scale=2, length=4, width=4]},
                 very thick,blue!40!black] (2.532,2.79) -- (2,1.08);         
  \draw (2.6, 1.9) node[fill=white] {$w$};                           
  
  \filldraw [gray] (2.532,2.79) circle (2pt);                        
  \draw (2.75,3.11) node[fill=white] {$x+v$};                        
  \draw [-{>[scale=2, length=4, width=4]},
                   very thick,green!40!black] (8,1.08) -- (2.532,2.79);                   
  \draw (3.7,2.7) node[fill=white] {$v$};                            
  
    \draw (0,3.3) node {(a)};
 
  
\end{tikzpicture}
\hfill
\begin{tikzpicture}[scale=0.68]
  \draw [ultra thick] (0,0) .. controls (2,2)  and (8,2) .. (10,0);  
  \draw (5,1.23) node[fill=white] {$M$};                             
  
  \filldraw [gray] (8,1.08) circle (2pt);                            
  \draw (8,.8) node[fill=white] {$x$};                               
  \draw [thick] (10,.47) -- (1,3.27);                                
  \draw (9.4,1.05) node[fill=white] {$T_x$};                         
  
  \filldraw [gray] (2,1.08) circle (2pt);                            
  \draw (2.8,.73) node[fill=white] {$y=x+v+w$};                      
  
  \draw [-{>[scale=2, length=4, width=4]},
                   very thick, blue!40!black] (2.532,2.79) -- (2,1.08);         
  \draw (2.6, 1.9) node[fill=white] {$w$};                           
  
  \filldraw [gray] (2.532,2.79) circle (2pt);                        
  \draw (2.75,3.11) node[fill=white] {$x+v$};                        
  \draw [-{>[scale=2, length=4, width=4]},
                    very thick,green!40!black] (8,1.08) -- (2.532,2.79);      
  \draw (3.7,2.7) node[fill=white] {$v$};                            
 
  \draw [thin] (0,.49) -- (10,3.46);                                 
  \draw (.4,1.05) node[fill=white] {$T_y$};                          
  \draw [-{>[scale=2, length=4, width=4]},
                   very thick,blue,dashed] (7.49959,2.71) -- (8,1.08);    
  \filldraw [gray] (7.49959,2.71) circle (2pt);                      
  \draw (7.35,3.11) node[fill=white] {$y+v^{\prime}$};               
  \draw [-{>[scale=2, length=4, width=4]},
                   very thick, green!80!black,dashed] (2,1.08) -- (7.49959,2.71);                             
  \draw (6.3,2.7) node[fill=white] {$v^{\prime}$};                   
  \draw (8.3, 1.9) node[fill=white] {$w^{\prime}$};                  
  
    \draw (0,3.3) node {(b)};
    
\end{tikzpicture}

\end{center}
\caption{Illustration of the sampler.
  (a) The vector $v \in T_x$ is projected normal to $T_x$ to $y = x + v + w\in M$.
   (b) The reverse step, needed for detailed balance, has $v^{\prime} \in T_y$, $w^{\prime} \in T_y^{\perp}$, 
   and $x = y + v^{\prime} + w^{\prime}$.
   The angle between $w$ and $T_y$ is the same as the angle between $w^{\prime}$ and $T_x$.
   }
\label{fig:prop}
\end{figure}

\subsection{The algorithm}

We describe the MCMC algorithm we use to sample a probability distribution of the form
\begin{equation}  \label{eq:pd}
   \rho(dx) = \frac{1}{Z} f(x)\, \sigma(dx) \; .
\end{equation}
Here $\sigma$ is the $d-$dimensional surface measure on $M$ given by \eqref{eq:manifold}.
As usual, the algorithm requires repeated evaluation of $f$ but does not use the unknown 
normalization constant $Z$.
The algorithm is described in pseudocode in Figure \ref{fig:MCMC}.

Our Metropolis style MCMC algorithm generates a sequence $X_k \in M$ with the property that 
\begin{equation}  \label{eq:bal}
   \mbox{if}   \;\;\;\;\;\;\;\; X_k \sim \rho \;\;\;\;\;\;\;\; 
   \mbox{then} \;\;\;\;\;\;\;\; X_{k+1} \sim \rho \; .
\end{equation}
Once $X_k = x$ is known, the algorithm makes a proposal $y \in M$.
If the proposal fails (see below) or if it is rejected, then $X_{k+1} = X_k$.
If the algorithm succeeds in generating a $y$, and if $y$ is accepted, then $X_{k+1} = y$.
If $M$ is connected, compact, and smooth, then the algorithm is geometrically ergodic.
See e.g. \cite{liu} for an efficient summary of the relevant MCMC theory.

The proposal process begins with a tangential move $x \to x+v$ with $v\in T_x$.
We generate $v$ by sampling a proposal density $p(v|x)$ defined on $T_x$.
All the numerical experiments reported here use an isotropic $d-$dimensional Gaussian with width
$s$ centered at $x$:
\begin{equation}\label{eq:pv}
        p(v|x) = \frac{1}{(2\pi)^{d/2}\,s^d}\,e^{-\frac{|v|^2}{2s^2}} \; .
\end{equation}
We generate $v$ using an orthonormal basis for $T_x$, which is the orthogonal
complement of the columns of the constraint gradient matrix $Q_x$, see (\ref{eq:Qx}).
This orthonormal basis is found as the last $d$ columns of the $d_a\times d_a$ matrix $Q$ in the QR decomposition\footnote{Please be aware of the
   conflict of notation here.
   The $d_a\times d_a$ orthogonal $Q$ of the QR decomposition is not the $d_a\times d$ gradient $Q_x$.} 
of $Q_x$, which we find using dense linear algebra.

Given $x$ and $v$, the projection step looks for $w \in T_x^{\perp}$ with $y = x + v + w \in M$. 
See Figure \ref{fig:prop}(a) for an illustration. 
It does this using an $m-$component column vector $a = (a_1, \ldots, a_m)^t$ and
\begin{equation}\label{eq:wx}
    w = \sum_{j=1}^m a_j \nabla q_j(x) = Q_x a \; .
\end{equation}
The unknown coefficients $a_j$ are found by solving the nonlinear equations
\begin{equation}\label{eq:proj}
  q_i\left(x+v+Q_xa\right) =0 , \qquad i=1,\ldots m.
\end{equation}
This can be done using any nonlinear equation solver. In our code, 
we solve \eqref{eq:proj} using simple Newton's method (no line search, regularization, etc.)
and initial guess $a = 0$.
The necessary Jacobian entries are
\begin{equation}  \label{eq:Ja}
       J_{ij} = \partial_{a_j} q_i(x + v + Q_xa) 
              = \left( \nabla q_i(x + v + Q_xa)\right)^t\,\nabla q_j(x) \; .
\end{equation}
We iterate until the first time a convergence criterion is satisfied
\begin{equation}\label{eq:converge}
|q(x+v+Q_xa)| \leq \mbox{\texttt{\textbf{tol}}} \; .
\end{equation}
where $|q|$ is some norm. In our implementation, we used the $l_2$ norm. 
When this criterion is satisfied, the projection is considered a success and $y = x + v + w$ is the proposal.
If the number of iterations reaches \texttt{\textbf{nmax}} before the convergence criterion is satisfied,
then the projection phase is considered a failure.
It is possible to change these details, but if you do you must also make the corresponding
changes in the detailed balance check below.

There are other possible projections from $x+v$ to $M$, each with advantages and disadvantages.
It would seem natural, for example, to find $y$ by computing
\begin{equation}  \label{eq:not}
       \argmin_{y \in M} \left| x + v - y \right|_2 \; .
\end{equation}
This produces a projection $w \in T_y^{\perp}$, rather than our $w \in T_x^{\perp}$.
It also has the advantage that such a $y$ always exists and is almost surely unique.
The optimization algorithm that computes (\ref{eq:not}) may be more robust than our 
Newton-based nonlinear equation solver.
The disadvantage of (\ref{eq:not}) is that there are curvature effects in the relative
surface area computation that are part of detailed balance. 
Computing these curvatures requires second derivatives of the constraint functions $q_i$. 
This is not be an intrinsic limitation, but requires more computation and is not as straightforward to implement as our algorithm. 

The proposal algorithm just described produces a random point $y \in M$.
There are two cases where we must reject it immediately so the detailed balance relation may hold (which we verify in section \ref{sec:db}.) 
First, we must check whether an inequality  constraint is violated, so $h_i(y) \leq 0$ for some $i$. If so, then $y$ is rejected.
Second, we must check whether the reverse step is possible, i.e. whether it is possible to propose $x$ starting from $y$. 
To make the reverse proposal, the algorithm would have to choose $v^{\prime} \in T_y$ so that 
$x = y + v^{\prime} + w^{\prime}$ with $w^{\prime} \perp T_y$.
Figure \ref{fig:prop}(b) illustrates this.
Since $x$ and $y$ are known, we may find $v^{\prime}$ and $w^{\prime}$ by the requirement that
$x-y = v^{\prime} + w^{\prime}$, $v^{\prime} \in T_y$, $w^{\prime} \in T_y^{\perp}$.
This is always possible by projecting $x-y$ onto $T_y$ and $T_y^\perp$, which are found using the $QR$ decomposition of $Q_y$.
However, we must additionally 
determine whether the Newton solver would find $x$ starting from $y + v^{\prime}$, a step that was neglected in earlier studies, e.g. \cite{lelievre2012}. 
Even though we know that $x$ is a solution to the equations, we 
 must still verify that the Newton solver produces $x$.
With that aim, we run the Newton algorithm from initial guess $y + v^{\prime}$ and see whether the 
the converge criterion \eqref{eq:converge} within $\mbox{\texttt{\textbf{nmax}}}$
iterations.
If not the proposal $y$ is rejected. Whenever $y$ is rejected we set $X_{k+1}=x$.

If we come to this point, we compute an acceptance probability $a(y|x)$ using the Metropolis Hastings formula
\begin{equation}\label{eq:a1}
     a(y|x) = \min \!\left( 1, \frac{f(y)p(v^{\prime}|y)}{f(x)p(v|x)}\right) \; .
\end{equation}
The reader might expect that we should use the more complicated formula
\[
     a(y|x) = \min \!\left( 1, \frac{f(y)p(v^{\prime}|y) J(x|y)}{f(x)p(v|x) J(y|x)}\right) \; ,
\]
where $J(y|x)$ is the inverse of the determinant of the linearization of the projection 
$v \rightarrow y$,  formally, $\left|\frac{\partial v}{\partial y} \right|$. This would account for how the area of a small patch of surface near $x$, 
is distorted upon mapping it to $y$. 
However, it turns out that the Jacobian factors cancel:
\begin{equation}  \label{eq:JJ}
        J(y|x)=J(x|y) \; .
\end{equation}
This is evident in part b of Figure \ref{fig:prop}, where the angle between $w$ and $T_y$ is the 
same as the angle between $w^{\prime}$ and $T_x$.
The proof in two dimensions is an exercise in Euclidean geometry.
In more dimensions, $J(y|x)$ and $J(x|y)$ are products of cosines of principle angles between
$T_x$ and $T_y$, see \cite{bjorck}.

To prove (\ref{eq:JJ}), let $U_x$ and $U_y$ be matrices whose columns are orthonormal bases 
for $T_x$ and $T_y$ respectively.
Column $i$ of $U_x$ will be denoted $u_{x,i}$, and similarly for $U_y$.
Consider a vector $\xi \in T_x$, which may be written $\xi = \sum a_i u_{x,i}$.
Let $\eta \in T_y$ be the projection of $\xi$ to $T_y$ normal to $T_x$, and write 
$\eta = \sum b_j u_{y,j}$.
The orthogonality condition is $u_{x,i}^t \eta = u_{x,i}^t\xi=a_i$.
This leads to the equations
\[
        \sum_j u_{x,i}^t u_{y,j}b_j = a_i \; .
\]
These take the form $Pb = a$, or $b = P^{-1}a$, where $P = U_x^t U_y$.
Since $U_x$ and $U_y$ are orthonormal bases, the volume element going from $x\to y$ is expanded by a factor
\begin{equation}  \label{eq:Jf}
      J(y|x)^{-1} = \mbox{det}\!\left(P^{-1}\right)= \mbox{det}\!\left((U_x^tU_y)^{-1}\right)\; .
\end{equation}
Similar reasoning shows that the projection from $T_y$ to $T_x$ perpendicular to $T_y$ 
expands volumes by a factor of
\begin{equation}
         J(x|y)^{-1} =  \mbox{det}\!\left((U_y^tU_x)^{-1}\right)\; .
\end{equation}
These factors are equal because the determinants of a matrix and its transpose are the same, proving  (\ref{eq:JJ}).

We close this section with a general remark concerning the choice of parameters. 
A key part of the algorithm is the projection step, which requires using Newton's 
method or some other algorithm to solve a nonlinear system of equations. 
This algorithm called twice at every step of the sampler and so making it efficient 
is critical to the overall efficiency of the sampler. 
Profiling reveals that most of the computing time is spent in the \texttt{\textbf{ project}} 
part (see Figure \ref{fig:MCMC}).
We found that setting the parameter \texttt{\textbf{nmax}}  low, or otherwise terminating the solver rapidly, made a significant difference to the overall efficiency. Unlike most problems that require solving equations, we are \emph{not} interested in guaranteeing that we can find a solution if one exists -- we only need to find certain solutions rapidly; our reverse projection check and subsequent rejection ensures that we correct for any asymmetries in the solver and still satisfy detailed balance.

\subsection{Verifying the balance relation}\label{sec:db}

In this section we verify what may already be clear, that the algorithm described above 
satisfies the desired balance relation (\ref{eq:bal}). 
 In the calculations below, we assume that the constraints hold \emph{exactly}, an assumption that will not be true when they are solved for numerically, as we discuss in the remark at the end of the section. 

To start, let us express the probability distribution of the proposal point $y$ as the product of a surface density and $d-$dimensional
surface measure, 
\begin{equation}
P_x(dy) = h(y|x) \sigma(dy) + \zeta(x) \delta_x(dy).
\end{equation}
Here $\delta_x(dy)$ is the unit point mass probability distribution at the point $x$, 
$\zeta(x)$ is the probability that the projection step failed, i.e. Newton's method failed to produce a $y$ or it lay outside the boundary, and $h(y|x)$ is the proposal density, i.e. the probability density of successfully generating $y\in M$ from $x$. 
The proposal density is the product of three factors:
\[
     h(y|x) = (1-{\bf 1}_{F_x}(y)) \;\cdot \; p(v|x) \;\cdot\; J(y|x) \; ,
\]
where the set $F_x\subset M$ consists of those points $y\in M$ that cannot be reached using our Newton solver starting from $x$,
and ${\bf 1}_{F_x}(y)$ denotes the characteristic function of the set $F_x$, which 
equals 1 if $y \in F_x$ and zero if $y \notin F_x$. 
Note that we must have $\int_M h(y|x) \,d\sigma(y)  + \zeta(x) = 1$.

If the acceptance probability is $a(y|x)$, then the overall probability distribution of $y$
takes the form 
\begin{equation}
       R_x(dy) = a(y|x) h(y|x) \sigma(dy) + \xi(x)\delta_x(dy) \; .
\end{equation}
Here, $\xi(x)$ is the overall probability of rejecting a move when starting at $x$, 
equal to the sum of the probability of not generating a proposal move, and the probability 
of successfully generating a proposal move that is subsequently rejected. 
Since $1-a(y|x)$ is the probability of rejecting a proposal at $y$, the overall
rejection probability is
\[
       \xi(x) = \zeta(x) + \int_M \left( 1 - a(y|x) \right)h(y|x)\sigma(dy) \; .
\]

We quickly verify that the necessary balance relation (\ref{eq:bal}) is satisfied if $a(y|x)$ is chosen to satisfy the detailed balance formula
\begin{equation}  \label{eq:db}
      f(x) a(y|x)h(y|x) = f(y) a(x|y)h(x|y) \; .
\end{equation}
For this, suppose $X_k$ has probability distribution on $M$ with density $g_k(x)$ with respect 
to surface area measure:
\[
    X_k \sim \rho_k(dx) = g_k(x) \sigma(dx) \; .
\]
We check that if (\ref{eq:db}) is satisfied, and $g_k = f$, then $g_{k+1} = f$. 
We do this in the usual way, which is to write the integral expression for $g_{k+1}$ and
then simplify the result using (\ref{eq:db}).
To start, 
\begin{equation*}
      g_{k+1}(x) = \int_M g_k(y) a(x|y) h(x|y) d\sigma(y) \;+\; \xi(x)g_k(x) \; .
\end{equation*}
The integral on the right represents jumps to $x$ from $y \neq x$ on $M$.
The second term on the right represents proposals from $x$ that were unsuccessful or rejected.
If $g_k = f$, using the detailed balance relation (\ref{eq:db}) this becomes
\begin{equation}\label{eq:gk+1}
      g_{k+1}(x) = \int_M f(x) a(y|x) h(y|x) d\sigma(y) \;+\; \xi(x)f(x) \; .
\end{equation}
The integral on the right now represents the probability
of accepting a proposal from $x$:
\[
     \mbox{Pr}\!\left(\,\mbox{accept from $x$}\right) = \int_M a(y|x)h(y|x) d\sigma(y) \; .
\]
When combined with the second term, the right side of \eqref{eq:gk+1} is equal to $f(x)$, as we wished to show.

We use the standard Metropolis Hastings formula to enforce (\ref{eq:db}):
\begin{equation}  \label{eq:r}
   a(y|x) = \min\left(1, \frac{f(y)h(x|y)}{f(x)h(y|x)}  \right) \; .
\end{equation}
If we succeed in proposing $y$ from $x$, then we know $y \notin F_x$ so we are able to evaluate $h(y|x)$. 
However, it is still possible that $x \notin F_y$, in which case $h(x|y) = 0$.
This is why we must apply the Newton solver to the reverse move, to see whether it
would succeed in proposing $x$ from $y$.
If not, we know $h(x|y) = 0$ and therefore $a(y|x) = 0$ and we must reject $y$.

\bigskip

\emph{Remark.} We have assumed so far the constraints hold exactly, which is not true when we solve for $y$ numerically. Rather, the proposal move will lie within some tolerance $\texttt{tol}$, so we are sampling a ``fattened'' region $M_{\texttt{tol}} = \left\{y:\;|q_i(y)| \leq \texttt{tol} \; i=1,\ldots,m,\; h_j(x) > 0, \; j = 1, \ldots, l \right\}$.  It is well-known (e.g. \cite{Fixman:1974dd,Hinch:1994ca}) that if the fattened region is uniformly sampled, the distribution near $M$ will not be the uniform distribution, but rather will differ by an $O(1)$ amount which does not vanish as $\texttt{tol}\to 0$. 

We do not know the distribution by which our numerical solver produces a point in $M_{\texttt{tol}}$. However, since we observe the correct distributions on $M$ in the examples below (for small enough \texttt{tol}), we infer that we are \emph{not} sampling the fattened region uniformly, but rather in a way that is consistent with the surface measure on the manifold, in the limit when $\texttt{tol} \to 0$.
Moreover, it is practical to take $\texttt{tol}$ to be extremely small, on the order of machine
precision, because of the fast local convergence of Newton's method.

\subsection{Examples}\label{sec:examples} 

In this section we illustrate the MCMC sampler with three examples: a torus, a cone and the special orthogonal group $SO(n)$.

\subsubsection{Torus}  \label{sec:torus}
Consider a torus $\mathbb{T}^2$ embedded in $\R^3$, implicitly defined by
\begin{equation}
\mathbb{T}^2 = \left\{ (x,y,z) \in \R^3 : \left( R- \sqrt{x^2+y^2}\right)^2 + z^2-r^2 = 0 \right\},
\end{equation}
where $R$ and $r$ are real positive numbers, with $R > r$. 
Geometrically, $\mathbb{T}^2$ is the set of points at distance
$r$ from the circle of radius $R$ in the $(x,y)$ plane centered at $(0,0,0)$.
An explicit parameterization of $\mathbb{T}^2$ is given by
\begin{equation}
\mathbb{T}^2 = \left\{ \left(R + r \text{cos}(\phi) \text{cos}(\theta), (R + r \text{cos}(\phi) \text{sin}(\theta), r \text{sin}(\phi) \right) : \theta, \phi \in [0,2\pi] \right\}.
\end{equation}

\begin{figure}[!t]
\begin{center}
\includegraphics[trim={10cm 2cm 6cm 2cm},clip,width=0.8\linewidth]{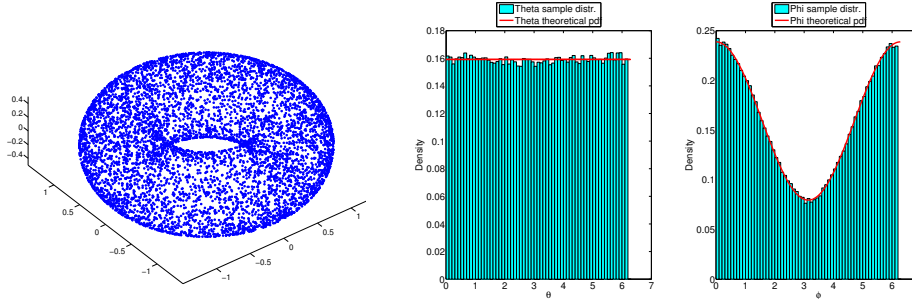}
\caption{Left: Points generated by our MCMC sampler on a torus. We generated $10^6$ points total and plot $10^4$ of these, one every $100$ MCMC steps.
Center and right: empirical and theoretical distribution of the torus angles $\theta$ (center) and $\phi$ (right).
The theoretical $\phi$ density (see, e.g.\ \cite{diaconis}) is 
$\frac{1}{2\pi}\left(1+\frac{r}{R} \cos(\phi) \right)$.}
\label{fig:torus_distr}
\end{center}
\end{figure}

We ran $N = 10^6$ MCMC steps to sample uniform measure ($f(x) = 1$) on $\mathbb{T}^2$
with toroidal radius $R = 1$, poloidal radius $r = .5$, and step size scale $s = .5$. 
Figure \ref{fig:torus_distr} shows our empirical marginal distributions and the exact theoretical distributions. 
The empirical distributions of the toroidal angle $\theta$ and the poloidal angle $\phi$
are correct to within statistical sampling error.
Around $6 \%$ of proposed moves were rejected because of failure in the reverse projection
step. When we didn't implement this step, the marginals were not correct. 

\subsubsection{Cone}\label{sec:cone}
 We consider the circular right-angle cone $\mathcal{C} \subset \R^3$
with vertex $(0,0,0)$ given by
\begin{equation}\label{cone}
    \mathcal{C} = \{ (x,y,z) \in \R^3 : z-\sqrt{x^2+y^2} = 0 , \; x^2+y^2 < 1, \; z > 0\}.
\end{equation}
We study this example in order to determine how our algorithm might behave near singular points 
of a hypersurface. 
We have written strict inequalities to be consistent with the general formalism (\ref{eq:manifold}).
In principle, the algorithm should do the same thing with non-strict inequalities, since the set of
points with equality (e.g. $z = 0$) has probability zero.
If we allowed $z=0$, then $\mathcal{C}$ would not be a differentiable manifold because of 
the singularity at the vertex.

\begin{figure}[!t]
\begin{center}
\includegraphics[trim = 160 40 10 10, clip, scale = 0.15]{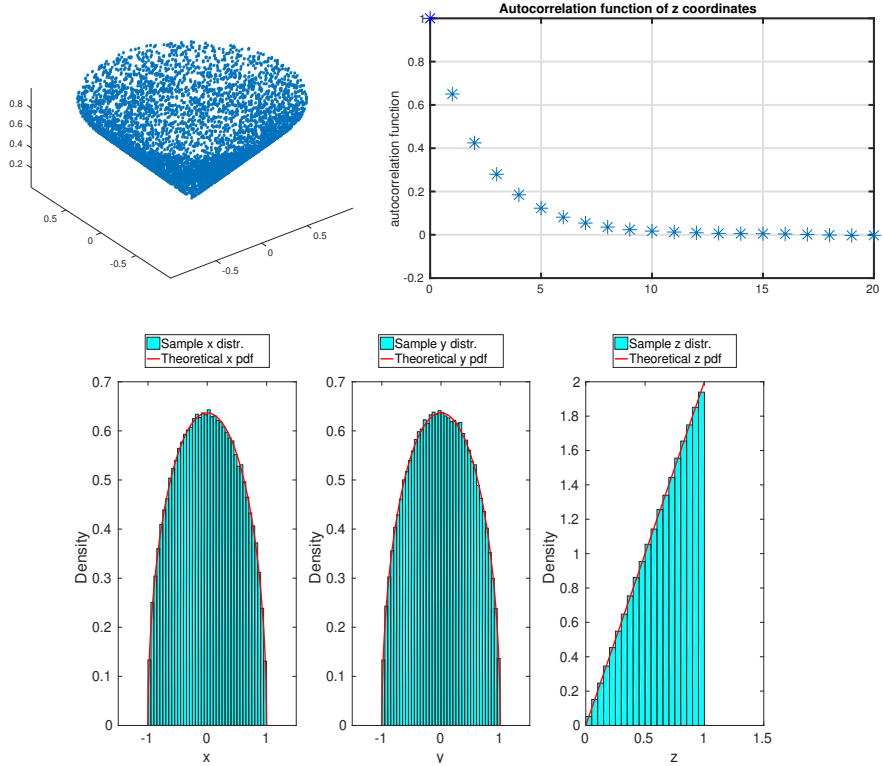}
\caption{An MCMC run of length $10^6$ for uniform sampling on a cone.  Top left: plot of $10^4$ points of the MCMC chain, over every $100$ steps. Top right: autocorrelation function of the coordinates $Z$.
Bottom: empirical and theoretical distributions for the coordinates $(X,Y,Z)$. 
}
\label{fig:cone_distr}
\end{center}
\end{figure}

We suspected that the projection steps (see Figure \ref{fig:MCMC}) would be failure prone near the 
vertex because of high curvature, and that the auto-correlation times would be large or even unbounded.
Nevertheless, Figure \ref{fig:cone_distr} shows satisfactory results, at least for this example.
We ran for $N=10^6$ MCMC steps with step size parameter $s = .9$. 
The theoretical marginal densities for $X$, $Y$, and $Z$ are easily found to be
\[
      g_X(x) = \frac{2}{\pi}\sqrt{1-x^2}\;,\;\;
      g_Y(y) = \frac{2}{\pi}\sqrt{1-y^2}\;,\;\;
      g_Z(z) =   2z \; .
\]

\subsubsection{Special orthogonal group} \label{sec:SO(n)}
We apply our surface sampler to the manifold $SO(n)$, which is the
the group of $n \times n$ orthogonal matrices with determinant equal to 1. 
We chose $SO(n)$ because these manifolds have high dimension, $\frac{1}{2}n(n-1)$, and high co-dimension, $\frac{1}{2}n(n+1)$ \cite{lee}. 
We view $SO(n)$ as the set of $n \times n$ matrices, $x \in \R^{n\times n}$ that satisfy the row 
ortho-normality constraints for $k = 1, \ldots, n$ and $l > k$:
\[
   g_{kk}(x) = \sum_{m=1}^n x_{km}^2 = 1 \;\; , \;\;\;\; g_{kl}(x) = \sum_{m=1}^n x_{km}x_{lm} = 0 \; .
\]
This is a set of $\frac{1}{2}n(n+1)$ equality constraints.
The gradient matrix $Q_x$ defined by (\ref{eq:Qx}) may be seen to have full rank for all $x \in SO(n)$.
Any $x$ satisfying these constraints has $\mbox{det}(x) = \pm 1$.
The set with $\mbox{det}(x) = 1$ is connected.
It is possible that our sampler would propose an $x$ with $\mbox{det}(x) = -1$, but we reject 
such proposals.
We again chose density $f(x) = 1$, which gives $SO(n)$ a measure that is, up to a constant, 
the same as the natural Haar measure on $SO(n)$.

One check of the correctness of the samples is the known large-$n$ distribution of $T = \mbox{Tr}(x)$,
which converges to a standard normal as $n \to \infty$ \cite{diaconisOn}.
This is surprising in view of the fact that $\mbox{Tr}(x)$ is the sum of $n$ random numbers each of
which has $\mbox{var}(x_{kk}) = 1$.
If the diagonal entries were independent with this variance, then $\mbox{var}(T) = n^{1/2}$
instead of the correct $\mbox{var}(T) = 1$.
This makes the correctness of the $T$ distribution an interesting test of the correctness 
of the sampler.

\begin{figure}[!t]
\begin{center}
\includegraphics[scale = 0.2]{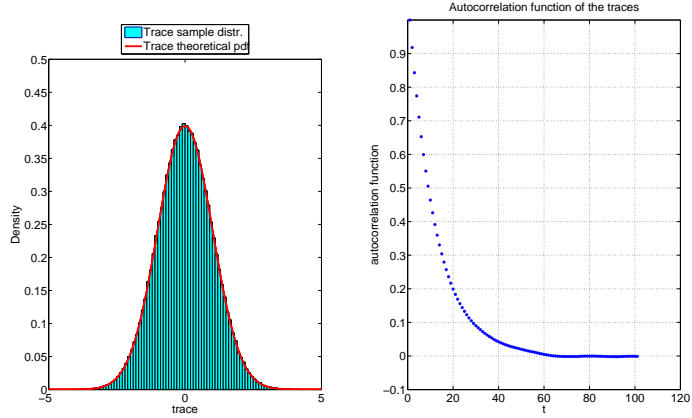}
\caption{An MCMC run of length $10^6$ for uniform sampling from $SO(11)$.
Left: the theoretical and empirical distribution of $T = \mbox{Tr}(x)$.
Right, the empirical auto-correlation of $T$.}
\label{fig:trace_distr}
\end{center}
\end{figure}

Figure \ref{fig:trace_distr} presents results for $n = 11$ and  $N=10^6$ MCMC steps, 
and proposal length scale $s = .28$. 
The distribution of $T$ seems correct but there are non-trivial auto-correlations.
The acceptance probability was about $35\%$.

Note that there are better ways to sample $SO(n)$ that to use our sampler.
For example, one can produce independent samples by choosing $y \in \R^{n\times n}$ with i.i.d. entries from ${\cal N}(0,1)$ and use the 
QR decomposition $y = Rx$, where $R$ is upper triangular and $x$ is orthogonal.

\section{Integrating over a manifold}  \label{sec:vol}

\subsection{Algorithm}\label{sec:volalg}

In this section we describe an algorithm to calculate integrals over a manifold $M$ of the form (\ref{eq:Z}).
We assume that $M$ given in \eqref{eq:manifold} is bounded. 
The strategy as outlined in the Introduction is to 
 consider a sequence of $d_a$-dimensional balls $B_0\supset B_1\supset \ldots \supset B_k$, centered at 
$x_0 \in M$ with radii $r_i$, for $i = 1,\ldots, k$, where $B_0\supset M$ entirely contains the manifold, 
and $B_k$ is ``small enough" in a way we describe momentarily. 
We \emph{assume} that each $M\cap B_i$ is connected, an assumption that is critical for our algorithm to work. 
We define the collection of integrals (see \eqref{eq:ZB})
\[
       Z_i = \int_{B_i \cap M} f(x) d\sigma(x) \; .
\]
The final estimate is $Z \approx \widehat{Z}$ according to (\ref{eq:Zhat}).

To estimate the ratio $R_i = Z_{i}/Z_{i+1}$ we generate $n_i$ points in $B_{i} \cap M$ with the probability distribution $\rho_i$ as in (\ref{eq:Jf}). 
This is done using the MCMC algorithm given in Section \ref{sec:MCMC} by including the extra inequality 
$|x-x_0|^2 < r_i^2$. 
An estimator for $R_i$ is then given by
\begin{equation}\label{eq:estimator_rho}
   \widehat{R}_i  = \frac{n_i}{N_{i,i+1}},
\end{equation}
where, for $i\leq j$, we define
\begin{equation}\label{nij}
N_{i,j} = \# \; \text{points generated in $B_i \cap M$ that also lie in $B_j \cap M$}.
\end{equation}
We could use additional elements $N_{i,j}$ with $j>i+1$ to enhance the estimate of the ratios, but we leave exactly how to do this a question for future research. Note that $N_{i,i}=n_i$.

\begin{figure}[!t]
\begin{center}
\includegraphics[scale = 0.12]{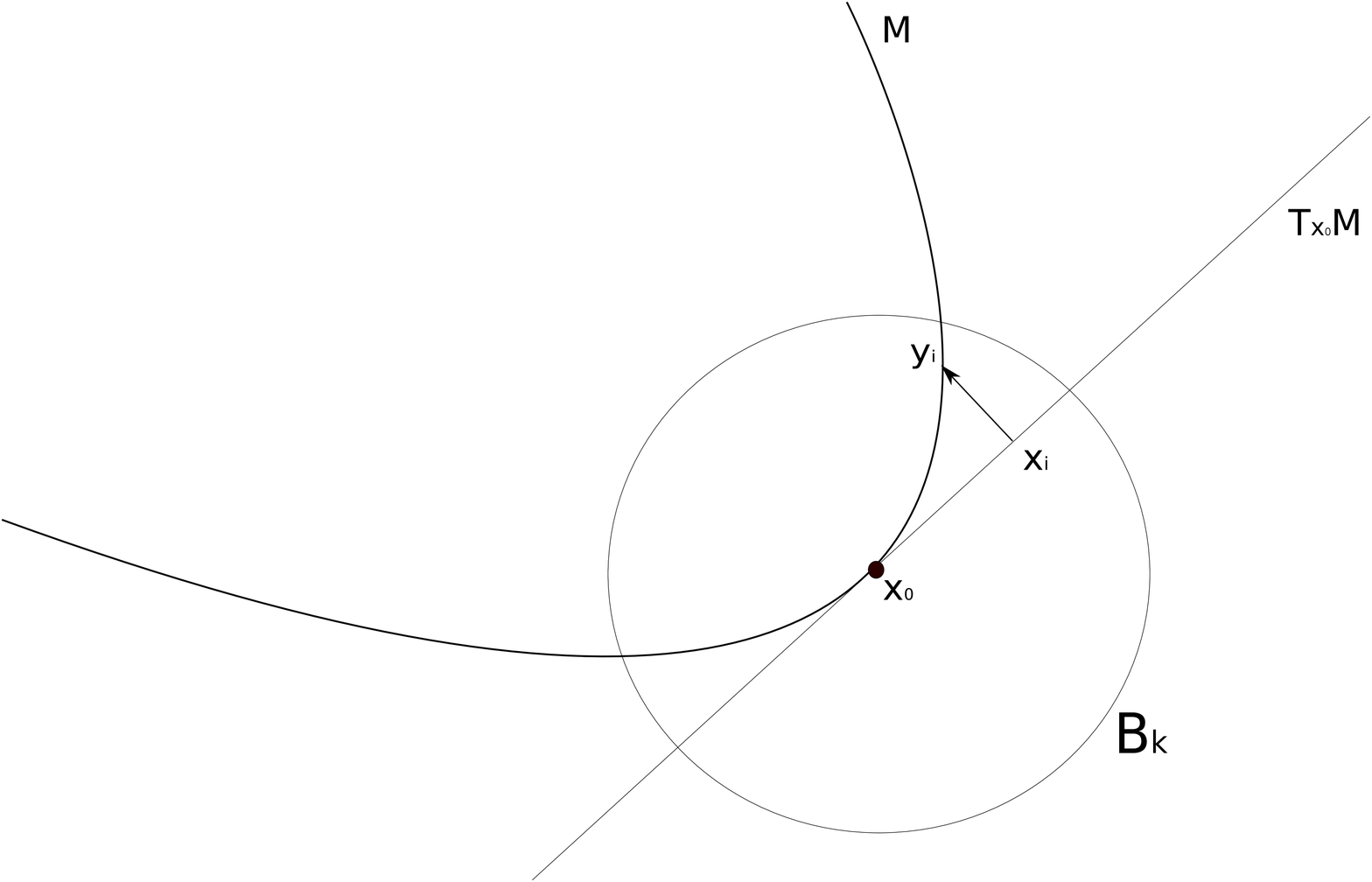}
\caption{Illustration of the projection step for the integration over the smallest intersection $B_k \cap M$ for a one-dimensional manifold: given a point $x_i$ uniformly generated in the ball $B_k \cap T_{x_0}M$, we find the unique point $y_i \in M$ by moving perpendicular to $T_{x_0}M$, starting from $x_i$. }
\label{fig:volume}
\end{center}
\end{figure}

The integral over the smallest set $M\cap B_k$ is computed in a different way.
If $M$ had no equality constraints, so that the dimension of $M$ were equal to the dimension of $B_k$, then we could choose 
$B_k\subset M$ to lie entirely inside $M$ and so the integral (for constant $f$) would be known analytically. 
Indeed, this idea was used in several studies that developed efficient algorithms to calculate 
the volume of a convex body \cite{simonovits, vempala}.  
With the equality constraints, in general there is no $B_k$ with radius $r_k > 0$ for which the integral $Z_k$ is known
in closed form exactly.

Instead, we choose $r_k$ small enough that 
that there is an easily computed one-to-one correspondence
between points in $M \cap B_k$ and points in the tangent plane $D_k = T_{x_0}\cap B_k$.
This happens, for example, when the manifold $M\cap B_k$ is a graph above $D_k$. 
We use the coordinates in the flat disk $D_k$ to integrate over $M\cap B_k$, which is simply the surface integral as illustrated 
in Figure \ref{fig:volume}:
\begin{equation}   \label{eq:ZD}  
    Z_k =  \int_{M\cap B_k} f(y) \,d\sigma(y) = \int_{D_k} f(y(x)) J^{-1}(x) \,dx \; .
\end{equation}
This uses the following notation.
The integral on the right is over the $d$-dimensional disk.
For any $x \in D_k$, the projection to $M$ perpendicular to $T_{x_0}$ is called $y(x)\in M$. 
For small enough $r_k$, the projection is easily calculated using the projection
algorithm of Figure \ref{fig:MCMC}.
The Jacobian factor $J(x)  = \text{det}(U_{x_0}^tU_{y(x)})$ is known from \eqref{eq:Jf}, 
where as before $U_x$ is a matrix whose columns are an orthonormal basis of $T_x$.

The integral over $D_k$ may be computed deterministically if $d$ is small enough, but otherwise 
we estimate it by direct Monte Carlo integration.
We choose $n_k$ independent points $x_i \in D_k$, which is easy for a $d-$dimensional disk.
We compute the corresponding projections and Jacobian weight factors, $y_i = y(x_i)$ and
$J_i = J(x_i)$.
Some of the $y_i$ will be outside of $M \cap B_k$, so we need an extra indicator function
factor, ${\bf 1}_{B_k}(y_i)$, that is equal to 1 if $y_i \in M\cap B_k$.
The $d$ dimensional volume of $D_k$ is known, so 
 we have the direct Monte Carlo estimate
\begin{equation}   \label{eq:Zke}   
        \widehat{Z}_k = 
           \frac{\mbox{vol}_d(D_k)}{n_k} \,\sum_{i=1}^{n_k} {\bf 1}_{B_k}(y_i)f(y_i) J_i^{-1} \; ,
           \qquad 
            \mbox{vol}_d(D_k) = \frac{\pi^{d/2}}{\Gamma(\frac{d}{2}+1)} \, r_k^d \; .
\end{equation}
The estimators (\ref{eq:estimator_rho}) and (\ref{eq:Zke}) are substituted in (\ref{eq:Zhat}) to 
give our overall estimate of the integral $Z$.
 
There is one additional consideration which is important to make this work in practice. That is, the projection from $D_k$ to $M\cap B_k$ must \emph{never fail}. If it does, the estimate \eqref{eq:Zke} will not be convergent. Empirically we found the projection always worked when $r_k$ was small enough, which proved to be a more limiting criterion than the criterion that $M\cap B_k$ be a graph above $D_k$.

\subsection{Variance estimate}\label{sec:errors}

This section serves two purposes.
One is to present the method we use to estimate the variance, $\sigma^2_{\widehat{Z}}$, 
of the estimator (\ref{eq:Zhat}).
This is similar to the strategy of \cite{hou}, but adapted
to this setting.
This estimate is accurate only if the MCMC runs are long enough so that $Z_k$ and the
$R_i$ are estimated reasonably accurately.
Our computational tests below show that the estimator can be useful in practice.
The other purpose is to set up notation for some heuristics for choosing parameters, a question dealt with in more depth in Section \ref{sec:toy}.

We start with the simple estimator $\widehat{Z}_k$ in \eqref{eq:Zke} for the integral over the smallest ball.
It is convenient to write the estimator as the sum of the exact 
quantity and a statistical error, and then to work with relative error rather than absolute error.
We neglect bias, as we believe that error from bias is much smaller than statistical error in all our computations.
In the present case, we write 
\begin{equation} \label{eq:Zkh}
      \widehat{Z}_k = Z_k + \sigma_k \zeta_k 
                   = Z_k \left( 1 + \rho_k\zeta_k\right) \; ,
\end{equation}
where $\zeta_k$ is a zero mean random variable with $\mbox{var}(\zeta_k) = 1$,
$\sigma_k$ is the absolute standard deviation of $\widehat{Z}_k$, and  
$\rho_k=\sigma_k/Z_k$ is the relative standard deviation.
Since the samples in (\ref{eq:Zke}) are independent, the standard deviation satisfies
\begin{equation*}
\sigma_k^2 =  \frac{\mbox{vol}^2_d(D_k)}{n_k} \text{Var}(G(y)), \qquad G(y) =  {\bf 1}_{B_k}(y)f(y) J^{-1}(y). 
\end{equation*}
We estimate $\text{Var}(G(y))$ with the sample variance and 
 obtain, for the absolute and relative variances, 
\begin{equation}\label{eq:sk}
\sigma_k^2 \approx \frac{\mbox{vol}^2_d(D_k)}{n_k^2} \sum_{i=1}^{n_k} (G(y_i)-\overline{G})^2, \qquad 
\rho_k^2 \approx \frac{\sigma_k^2}{\widehat{Z}_k^2}. 
\end{equation}
Here $\overline{G} = \sum_{i=1}^{n_k}\frac{1}{n_k} G(y_i)$ is the sample mean. 

Now consider the estimator $\widehat{R}_i$ in \eqref{eq:estimator_rho} for the $i$th volume ratio. This involves the number count $N_{i,i+1}$. We start by approximating  $\mbox{var}(N_{i,i+1})$ and then incorporate this in the approximation for $\mbox{var}(\widehat{R}_i)$. 
The number count may be written
\begin{equation}   \label{eq:Ni}
         N_{i,i+1} = \sum_{j=1}^{n_i} {\bf 1}_{B_{i+1}}(X_j) \; ,
\end{equation}
where $X_j \in M \cap B_i$ is our sequence of samples. 
In the invariant probability distribution of the chain, the probability of $X_j \in B_{i+1}$ is
\begin{equation}  \label{eq:pi}
        p_i = \mathbb{E}\!\left[ {\bf 1}_{B_{i+1}}(X_j)\right] = \frac{1}{R_i} = \frac{Z_{i+1}}{Z_i} \; .
\end{equation}
Since ${\bf 1}_{B_{i+1}}(X_j)$ is a Bernoulli random variable, we have (in the invariant distribution)
\begin{equation}  \label{eq:vpi}
   \mbox{var}\!\left( {\bf 1}_{B_{i+1}}(X_j)\right) = p_i (1-p_i) \; .
\end{equation}
We will call this the static variance. It alone is not enough to determine $\mbox{var}(N_{i,i+1})$, since the steps of the Markov chain are correlated. To account for this correlation, 
we use the standard theory (see e.g.\ \cite{sokal}) of error bars for MCMC estimators of 
the form (\ref{eq:Ni}), which we now briefly review. 

Consider a general sum over a function of a general MCMC process
\[
        S = \sum_{j=1}^n F(X_j) \; .
\]
The equilibrium lag $t$ auto-covariance is (assuming $X_j$ is in the invariant distribution
and $t \geq 0$)
\[
      C_t = \mbox{cov}\!\left[\,F(X_j),F(X_{j+t})\right] \; .
\]
This is extended to be a symmetric function of $t$ using $C_t = C_{\left|t\right|}$.
It is estimated from the MCMC run using
\begin{equation}  \label{eq:Ct}
    \overline{F} \approx \frac{1}{n_i} \sum_{j=1}^{n_i} F(X_j) \; , \;\;\;\;
    C_t \approx  \frac{1}{n_i-t} \sum_{j=1}^{n_i-t} 
          \left( F(X_j)- \overline{F}\right)\left( F(X_{j+t})- \overline{F}\right) \; .
\end{equation}
The Einstein Kubo sum is
\[
         D = \sum_{t = - \infty}^{\infty} C_t \; .
\]
The auto-correlation time is
\begin{equation}  \label{eq:tau}
        \tau = \frac{D}{C_0} = 1 + \frac{2}{\mbox{var}(F(X_j))} \sum_{t=1}^{\infty} C_t \; .
\end{equation}
The sum may be estimated from data using the self-consistent window approximation described in \cite{sokal}.
The static variance $C_0$ in (\ref{eq:tau}) may be estimated using 
(\ref{eq:vpi}) rather than (\ref{eq:Ct}).
The point of all this is the large $n$ variance approximation
\[
   \mbox{var}(S) \approx nD = nC_0\tau \; .
\]

For our particular application, this specializes to
\begin{equation}  \label{eq:vNi}
   \mbox{var}(N_{i,i+1}) \approx  np_i(1-p_i)\tau_i \; ,
\end{equation}
where $\tau_i$ is the correlation time for the indicator function ${\bf 1}_{B_{i+1}}(X_j)$, assuming points are generated in $M\cap B_i$. It may be calculated using \eqref{eq:tau}. 

We use this to write $N_{i,i+1}$ in terms of a random variable $\xi_i$ 
with unit variance and mean zero as 
\begin{align}
      N_{i,i+1} 
                \approx n_ip_i \left( 1 + \sqrt{ \frac{(1-p_i)\tau_i}{n_ip_i}}\,\xi_i \right) \; .
    \label{eq:Nierr}
\end{align}
Substituting this approximation into \eqref{eq:estimator_rho} gives the approximation
\begin{equation}\label{eq:Rapprox}
      \widehat{R_i}    \approx \frac{1}{p_i} \left( 1 -   \sqrt{ \frac{(1-p_i)\tau_i}{n_ip_i}}\,\xi_i \right) \; .
\end{equation}
where we expanded the term $1/N_{i,i+1}$ assuming that $n_i$ was large. 
Continuing, we substitute into (\ref{eq:Zhat}) to get
\begin{align*}
     \widehat{Z} 
                 \approx Z \left[\left( 1 + \rho_k \zeta_k\right)
                     \prod_{i=0}^{k-1} \left( 1 -   \sqrt{ \frac{(1-p_i)\tau_i}{n_ip_i}}\,\xi_i \right)\right]\;.
\end{align*}
If the relative errors are all small, then products of them should be smaller still.
This suggests that we get a good approximation of the overall relative error keeping only
terms linear in $\zeta_k$ and the $\xi_i$.
This gives
\[
      \widehat{Z} \approx Z 
         \left[ 1 + \rho_k \zeta_k + \sum_{i=0}^{k-1} \sqrt{ \frac{(1-p_i)\tau_i}{n_ip_i}}\,\xi_i \right] \; .
\]
Finally, we assume that the random variables $\zeta_k$ and $\xi_i$ are all independent.
As with the earlier approximations, this becomes valid in the limit of large $n_i$.
In that case, we get
\begin{equation}\label{eq:varZ}
    \mbox{var}(\widehat{Z}) 
      \;=\; \sigma_{\widehat{Z}}^2   \;\;\; \approx\;\;\;  Z^2 \left( \rho_k^2 + \sum_i \frac{(1-p_i)\tau_i}{n_ip_i} \right)  \; .
\end{equation}
The one standard deviation error bar for $\widehat{Z}$ is $Z\sigma_r$, where the relative 
standard deviation is
\begin{equation}    \label{eq:sr}  
      \sigma_r \approx \sqrt{\rho_k^2 + \sum_{i=0}^{k-1} \frac{(1-p_i)\tau_i}{n_ip_i}} \; .
\end{equation}

All the quantities on the right in \eqref{eq:sr} are estimated from a single run of the sampler in each $M\cap B_i$, with $p_i$ estimated as $1/\widehat{R}_i$  (see \eqref{eq:estimator_rho}), $\tau_i$ estimated from 
the sequence  ${\bf 1}_{B_{i+1}}(X_j)$ using a self-consistent window approximation to (\ref{eq:tau}), 
and $\rho_k$ estimated as in \eqref{eq:sk}.

\subsection{Algorithm parameters}\label{sec:implement}

This section describes how we choose   
the base point, $x_0$, and the sequence of ball radii $r_0, r_1,\ldots,r_k$. 
A guiding motivation is the desire to minimize the variance of $\widehat{Z}$, 
estimated by \eqref{eq:varZ}, for a given amount of work. 
Our parameter choices are heuristic and certainly may be improved.
But they lead to an overall algorithm that is ``correct'' in that it gives convergent estimates 
of the integral $Z$.

Our aim is to choose $r_0$ as small as possible, and $r_k$ as large as possible, to have a smaller overall ratio to estimate. Additionally, we expect the estimate of the final integral $\widehat{Z}_k$ to be significantly faster and more accurate than each ratio estimate $\widehat{R}_i$, since $\widehat{Z}_k$ is estimated by independently chosen points, so we wish $r_k$ to be large. 
With this aim in mind, the first step is to learn about $M$ through sampling.
We generate $n$ points $x_i \in M$ by sampling the surface measure 
$\rho(x) = \frac{1}{Z} f(x) d\sigma(x)$, using $n$ MCMC steps of the surface sampling algorithm of Section \ref{sec:MCMC}.
Then we choose $x_0$ from among these samples.  
Usually we would like it to be near the center of $M$, to reduce $r_0$ or increase $r_k$. Unless otherwise noted, for our implementation we choose it far from an inequality constraint boundary using 
\begin{equation}  \label{eq:O1}
   x_0 = \argmax_{x_i = x_1, \ldots, x_n}\;\; \min_j \; h_j(x_i) \; .
        \hspace{60pt}
\end{equation}
The $h_j$ are the functions that determine the boundary of $M$ (cf. \eqref{eq:manifold}).
This choice may allow the smallest ball, $B_k$ to be large, or make it so that most points of
$B_k$ project to points that satisfy the inequality constraints.
This is heuristic in that there is no quantitative relation between $\min_j h_j(x_i)$ and the distance (Euclidean or geodesic or whatever) from $x_i$ to the boundary of $M$, if $M$ has a 
boundary.

This is not the only possibility for choosing $x_0$. One could also choose it randomly from the set of sampled points. Another possibility would be to center $x_0$ by minimizing the ``radius'' of $M$ about $x_0$, which is
the maximum distance from $x_0$ to another $x \in M$. 
We could do this by setting $
   x_0 = \argmin_{x_i}
   \max_j \left| x_i - x_j\right| \; .
$
Having a smaller radius may mean that we need fewer balls $B_j$ in the integration 
algorithm.

Once $x_0$ is fixed, the radius $r_0$ of the biggest ball $B_0$ can be found using the sample points: 
\begin{equation}  \label{eq:r0}
   r_0 = \max \{ |x_0-x_i| : i = 1, \ldots, n \}.
\end{equation}



Next consider the minimum radius $r_k$. We need this to be small enough that (a) $B_k\cap M$ is 
single-valued over $T_{x_0}$,  and (b) the projection from $D_k\to M\cap B_k$ never fails. 
We test (b) first because it is the limiting factor. 
We choose a candidate radius $\tilde r$, sample $n\approx 10^5$ points uniformly from $D_k$, and project to $M\cap B_k$. If any single point fails, we shrink $\tilde r$ by a factor of 2 and test again. 
When we have a radius $\tilde r$ where all of the $10^5$ projections succeed, our error in estimating \eqref{eq:Zke} should be less than about 0.001\%. 

Then, we proceed to test (a) approximately as follows. 
We sample $n$ points $x_1, \ldots, x_n$ uniformly from $B_{\tilde{r}}(x_0) \cap M$, where $B_{\tilde r}(x_0)$ is the ball of radius $\tilde r$ centered at $x_0$. 
We then consider every pair $(x_i,x_j)$ and check the angle of the vector $v_{ij} = x_i-x_j$ with the orthogonal complement $T_{x_0}^{\perp}$. 
Specifically, we construct an orthogonal projection matrix onto the normal space $T_{x_0}^{\perp}$, for example as $P = U_{x_0}U_{x_0}^T$. 
Then, if there exists a pair for which $|v-Pv| < \texttt{tol}$, where $\texttt{tol}$ is a tolerance parameter, we shrink $\tilde r$ and repeat the test. In practice, we shrink $\tilde r$ by a constant factor each time (we used a factor $C = 2$). 
We expect this test to detect regions where $B_{\tilde{r}}(x_0) \cap M$ is multiple-valued as long as it is sampled densely enough. Of course, if the overlapping regions are small, we may miss them, but then these small regions will not contribute much to the integral $Z_k$. 

An alternative test (which we have not implemented) would be to sample $B_{\tilde{r}}(x_0) \cap M$ and check the sign of the determinant of the projection to $T_{x_0}$, namely $\mbox{sgn}\left( \mbox{det}(U_{x_0}^tU_y)\right)$ where $y\in B_{\tilde{r}}(x_0) \cap M$. If $B_{\tilde{r}}(x_0) \cap M$ is multiple-valued, we expect there to be open sets where this determinant is negative, and therefore we would find these sets by sampling densely enough. 

Having chosen $x_0$, $r_0$, and $r_k$, it remains to choose the intermediate ball radii $r_i$.
We use a fixed $\nu>1$ and take the radii so the balls have fixed $d$-dimensional volume ratios (note $d$ is typically different from the ball's actual dimension $d_a$), as
\begin{equation}      \label{eq:nu}
        \left(\frac{r_{i}}{r_{i+1}} \right)^d = \nu \; .
\end{equation}
Such fixed $\nu$ strategies are often used for volume computation \cite{vempala}. 
Given (\ref{eq:nu}), the radii are
$        r_i = r_0 \nu^{-i/d}  .$
Since we have already chosen the smallest radius, this gives a relationship between $\nu$
and $k$ as
        $ \nu = \left(r_0/r_k\right)^{d/k}$. 
Note that we choose $r_k$, the smallest radius, before we choose $k$, the number of stages.

It is not clear what would be the best choice for $\nu$.
If $\nu$ is too large then the ratios $R_i$ (see \eqref{eq:estimator_rho})
will not be estimated accurately -- indeed they will be large so $p_i$ in \eqref{eq:sr}  will be small, increasing the relative standard deviation of the estimator. 
If $\nu$ is too small, then $k$ will be large and the errors at each stage will accumulate. In addition, the error at each single stage will be large if we fix the total number of MCMC points, since there will be fewer points per stage.
Clearly, these issues are important for the quality of the algorithm
in practice, but they are issues that are seriously under-examined in the literature. 
We return to them in section \ref{sec:toy}. 

\section{Integration: examples}\label{sec:results} 

This section presents some experiments with the integration algorithm of Section \ref{sec:vol}.
There are three sets of computations.
We first compute the surface area of the torus $\mathbb{T}^2$ (section \ref{sec:torusvol}.)
This computation is cheap enough that we can repeat it many times and verify that the 
error bar formula (\ref{eq:sr}) is reasonably accurate.
We then calculate the volume of $SO(n)$ (section \ref{sec:SO(n)vol}), to verify 
the algorithm can work on surfaces with high dimension and high co-dimension.
Finally, we apply the algorithm to study clusters of sticky spheres (section \ref{sec:spheres}), and make new predictions that could be tested experimentally.  

\subsection{Surface of a torus}\label{sec:torusvol}
 Consider the torus of Subsection \ref{sec:torus} with $R=1$
and $r = .5$.
The exact area is $Z = 4 \pi^2 rR \approx 19.7392$.
The geometry of the manifold is simple enough that we can choose the parameters of the algorithm analytically. We take the initial point to be $x_0 = (R+r,0,0)$. 
The smallest ball $B_0$ centered at $x_0$ such that $\mathbb T^2\subset \mathbb T^2\cap B_0$ has radius $r_0 = 2(R+r)$.
We take the radius of the smallest ball to be $r_k = r$, for which the projection 
$T_{x_0}\cap B_k \to \mathbb{T}^2\cap B_k$ is clearly single valued.
For each computation we generate $n$ points at each stage, with the total number of points fixed
to be $n_{t}= kn = 10^5 $.

\begin{figure}[!t]
\begin{center}
\includegraphics[scale = 0.2]{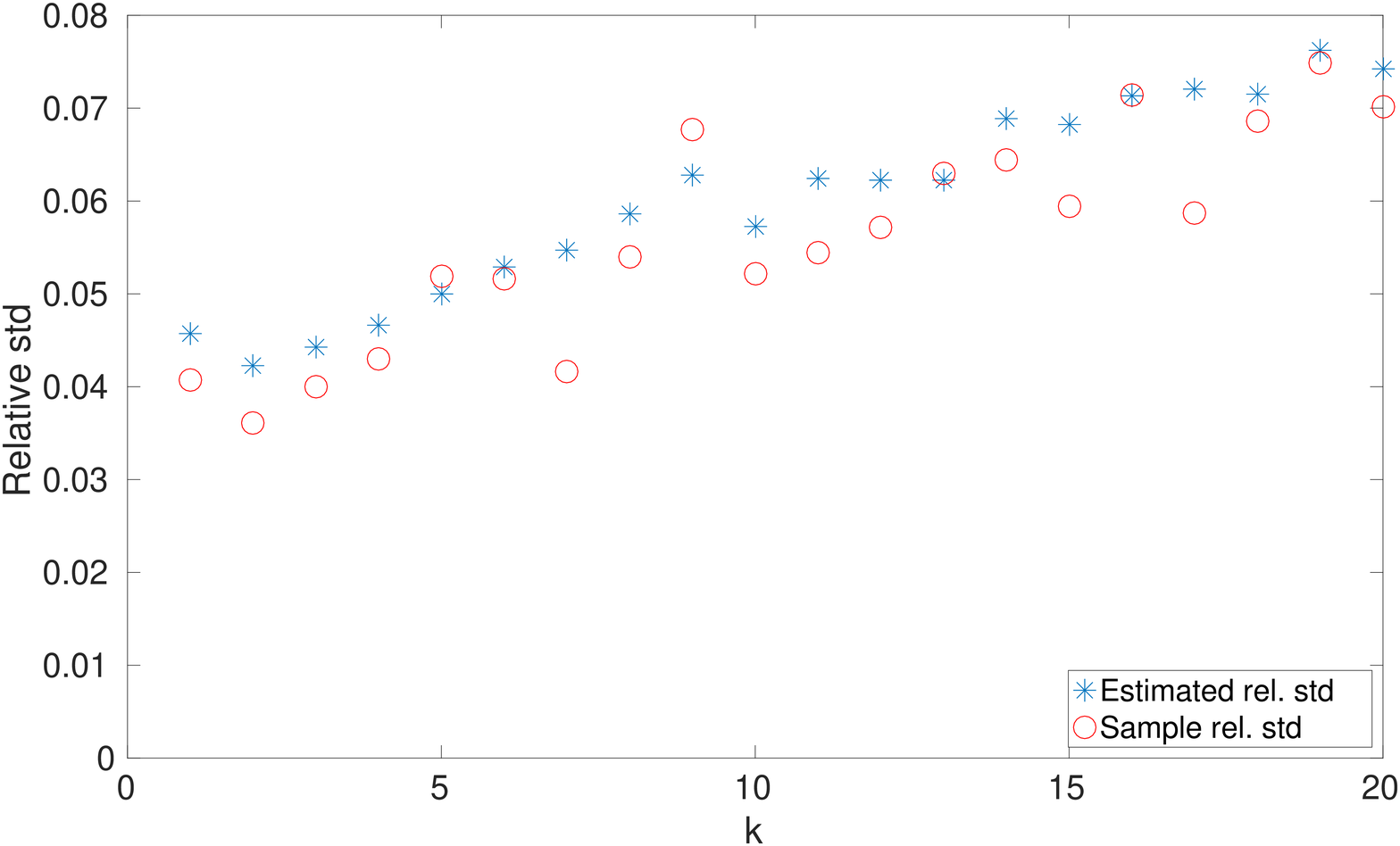}
\caption{Relative standard deviations $\sigma_r$ for the volume estimator on a torus, calculated using a total of $10^5$ MCMC points with other parameters given in section \ref{sec:torusvol}. The standard deviation was estimated using our formula \eqref{eq:sr} applied to one run (blue stars), and also 
as the sample standard deviation of $50$ independent runs (red circles). 
}
\label{fig:error_torus}
\end{center}
\end{figure}

We calculated the volume for each number of stages $k = 1, \ldots, 20$.
We repeated each calculation $50$ times independently.
This allowed us to get a direct estimate of the variance of the volume estimator $\widehat{Z}$.
Figure \ref{fig:error_torus} shows that the standard deviation estimated from a single run using (\ref{eq:sr})
is in reasonable agreement with the direct estimate from $50$ independent runs.\footnote{We do
      not strive for very precise agreement.
      An early mentor of one of us (JG), Malvin Kalos, has the wise advice: ``Don't put error bars 
      on error bars.''}
It also shows that for this problem $k=2$ seems to be optimal.
It does not help to use a lot of stages, but it is better to use more than one.

\subsection{The special orthogonal group} \label{sec:SO(n)vol}
The volume of $SO(n)$ as described in Subsection \ref{sec:SO(n)} is (see, e.g., \cite{volumes})
\begin{equation*}
\mbox{vol}(SO(n)) = 2^{\frac{n(n-1)}{4}} \prod_{i=1}^{n-1} \mbox{vol}(S^i),
\qquad \text{where}\quad
\mbox{vol}(S^i) = \frac{2\pi^{\frac{i+1}{2}}}{\Gamma\left(\frac{i+1}{2}\right)}
\end{equation*} 
is the volume of the $i$th dimensional sphere,. 
This formula may be understood as follows.
An element of $SO(n)$ is a positively oriented orthonormal basis of $\mathbb{R}^n$.
The first basis element may be any element of the unit sphere in $\mathbb{R}^n$, which is $S^{n-1}$.
The next basis element may be any element orthogonal to the first basis element, so it is in
$S^{n-2}$, and so on.
The prefactor $2^{\frac{n(n-1)}{4}}$ arises from the embedding of $SO(n)$ into 
$\mathbb{R}^{n\times n}$.
For example, $SO(2)$ is a circle that would seem to have length $2\pi$, but 
in $\mathbb{R}^4$ it is the set $(\cos(\theta), \sin(\theta), -\sin(\theta), \cos(\theta))$,
$0 \leq \theta < 2\pi$.
This curve has length $2 \sqrt{2}\pi$.
The difference is the ratio $2^{\frac{2\times 1}{4}} = \sqrt{2}$.

We applied the integration algorithm to compute the volume of $SO(n)$, for $2\leq n \leq 7$. 
For each $n$, we ran the algorithm $N = 50$ independent times and computed the mean 
of the  estimates $\overline{V}$, and the sample standard deviation $s_V$. 
We used a total of $n_t = 10^5$ points and $k = 4$ steps for each $n$.  Since the true value $V_t$ is known, we can compute the relative standard deviation $\widetilde{s_v}$ and the relative error $\epsilon_r$ via
\begin{equation}\label{eq:errors}
\widetilde{s_V} = \frac{s_V}{\sqrt{N} V_t}, \qquad \epsilon_r = \left|\frac{\overline{V}-V_t}{V_t}\right|. 
\end{equation}
We give the results in Table \ref{orth_group}, which shows our relative errors are small, and close the the estimated standard deviations. 
We notice that the errors increase as the dimension gets higher.  

\begin{table}[!t]
\begin{center}
\begin{tabular}{c c c c r@{.}l r@{.}l}
$n$ & $d$ & $V_t$ & $\overline{V}$  & \multicolumn{2}{c}{$\widetilde{s_V}$} & \multicolumn{2}{c}{$\epsilon_r$} \\
\hline
2 & 1 & 8.89 & 8.9 & 0&0038 & 0&0047 \\
3 & 3 &  223.3 & 223 & 0&0078 & 0&0024 \\
4 & 6 & $1.24 \cdot 10^4$ & $1.2 \cdot 10^4$ & 0&016 & 0&012 \\
5 & 10 & $1.31 \cdot 10^6$ & $1.3 \cdot 10^6$ & 0&022 & 0&024 \\
6 & 15 & $2.30 \cdot 10^8$ & $2.2 \cdot 10^8$ & 0&0354 & 0&0519 \\
7 & 21 & $6.09 \cdot 10^{10}$ & $6.0 \cdot 10^{10}$ & 0&0547 & 0&0094
\end{tabular}
\end{center}
\caption{Numerical vs analytical results for the volume of the special orthogonal group $SO(n)$, 
for $2 \leq n \leq 7$. Here $d$ is the dimension of $SO(n)$, $\overline{V}$ is the sample mean obtained from a sample of size $50$ and the relative standard deviation $\widetilde{s_V}$ and relative error $\epsilon_r$ are given in \eqref{eq:errors}.
}
\label{orth_group}
\end{table}

\subsection{Sticky-sphere clusters}\label{sec:spheres}

Next we consider a system of particles interacting with a so-called ``sticky'' potential, which can be thought of as a delta-function when the surfaces of two particles are in contact. This is a model for particles whose range of interaction is very short compared to the particles' diameters, as is the case for many kinds of colloidal particles. The probabilities of finding the system in different states in equilibrium can be calculated using statistical mechanics \cite{mhc2017} and we briefly summarize the relevant ideas. 

Consider a system of $N$ unit spheres in $\mathbb R^3$, represented as a vector 
$x = (x_1, \ldots, x_N)^t \in \R^{3N}$, where $x_i \in \R^3$ is the center of the $i$-th sphere. 
We suppose there are $m$ pairs of spheres in contact $E = \{(i_1,j_1), \ldots, (i_m,j_m)\}$. 
For each pair in contact there is a constraint
\begin{equation}\label{constraints}
     |x_i-x_j|^2-1 = 0, \qquad (i,j)\in E,
\end{equation} 
and for each pair not in contact there is an inequality which says the spheres cannot overlap:
\begin{equation}\label{boundary}
        |x_i - x_j|> 1, \qquad (i,j)\notin E. 
\end{equation}
To remove the translational degrees of freedom of the system, we fix the center of mass of the 
cluster at the origin.
We do this by imposing three extra constraints, which together are written in vector form as 
\begin{equation}\label{con2}
     \sum_{i=1}^N x_i = 0 \; .
\end{equation}

We define $M$ to be the set of points $x$ satisfying the constraints \eqref{constraints}, 
\eqref{con2}, and the inequalities \eqref{boundary}. 
It is possible that $M$ has singular points where the equality constraints \eqref{constraints} 
are degenerate (have gradients that are not linearly independent), but we ignore this possibility (note that the cone in section \ref{sec:cone} shows our algorithm may work even near singular points.) Therefore $M$ is a manifold of dimension $d = 3N - m -3$ embedded in an ambient space of dimension $d_a = 3N$, with surface area element $d\sigma$. 
If $M$ is not connected, our algorithm samples a single connected component. 
 
The equilibrium probability to find the spheres in a particular cluster defined by contacts $E$ is proportional to the partition function $Z$. In the sticky-sphere limit, and where the cluster is in a container much larger than the cluster itself, the partition function is calculated (up to constants that are the same for all clusters with the same number of spheres) as \cite{miranda,mhc2017}
\begin{equation}\label{eq:Zsticky}
   Z = \kappa^m z, \qquad \text{where} \quad z = \int_{M} f(x)\, d\sigma(x).
\end{equation}
Here $\kappa$ is the ``sticky parameter'', a single number depending on both the strength of the interaction between the spheres and the temperature. This parameter is something that must be known or measured for a given system. 
The factor $z$ is called the ``geometrical partition function'' because it does not depend on the nature of the interaction between the spheres, nor the temperature. 
The function $f(x)$ is given by
\begin{equation}\label{eq:f}
f(x) = \prod_{i=1}^{3N-m-3}\lambda_i(x)^{-1/2},
\end{equation}
which is the product of the non-zero eigenvalues $\lambda_i$ of $R^tR$, where $R$ is half the Jacobian of the set of contact constraints \eqref{constraints}. In the mathematical theory of rigidity, the matrix $R$ is called the rigidity matrix \cite{Connelly:1996vj}. 

As a first check that our algorithm is working we calculate geometrical partition functions $z$ that are already known, by considering the 13 different clusters of $6$ spheres with $m = 10$ contacts. All manifolds formed in this way are five-dimensional.  
For each manifold we compute $z$ as in \eqref{eq:Zsticky} using $n_t = 10^8$ points and $k = 4$ steps. For each estimate we compute its error bar using formula \eqref{eq:sr}.   We compare these values with the values $z_t$ calculated in \cite{miranda}, which were obtained by  parameterizing and triangulating each manifold
and calculating the integrals using finite elements.\footnote{
In \cite{miranda}, the calculations were done on the quotient manifold obtained by modding out the rotational degrees of freedom. We originally tried working in a quotient space here, by fixing six coordinates; this worked for these low-dimensional manifolds but we had problems with higher-dimensional ones, particularly chains. We presume this was because of the singularity in the quotient space. To be consistent, all calculations in this paper are performed on the original manifold $M$, which contains rotational degrees of freedom.
Therefore, we compare our results to the value $z_t = 8 \pi^2 \widetilde{z}$, where $\widetilde{z}$ is the value reported in \cite{miranda}. The constant $8\pi^2$ is the volume of the rotational space $SO(3)$.
} 
Table \ref{2dmodes} shows the two calculations agree.

\begin{table}[!t]
\begin{center}
\begin{tabular}{c r@{.}l @{\,} c@{\,} l   r@{.}l  }
Manifold &  \multicolumn{4}{c}{$z$}   &  \multicolumn{2}{c}{$z_t$}\\
\hline
$8$ & 6 & 99 &$\pm$& 0.12 & 7 & 014   \\
$9$ & 44 & 4 &$\pm$& 0.33 & 44 & 402 \\
$10$ & 57&1 &$\pm$& 0.26 & 57 & 07 \\
$11$ & 18&9 &$\pm$& 0.20 & 19 & 18 \\
$12$ & 3&30 &$\pm$& 0.020 & 3 & 299 \\
$13$ & 62&4 &$\pm$& 0.23 & 62 & 34 \\
$14$ & 3&35 &$\pm$& 0.039 & 3 & 299 \\
$15$ & 27&7 &$\pm$& 0.14 & 27 & 78 \\
$16$ & 11&7 &$\pm$& 0.058 & 11 & 67 \\
$17$ & 28&5 &$\pm$& 0.31 & 28 & 41 \\
$18$ & 8&44 &$\pm$& 0.088  & 8 & 238 \\
$19$ & 69&2 &$\pm$& 0.33 & 69 & 43 \\
$20$ & 50&9 &$\pm$& 0.27 & 50 & 86 \\
\end{tabular}
\end{center}
\caption{Geometrical partition functions $z$ of the configuration spaces of clusters of $N=6$ 
unit spheres with $m=10$ contacts (the numbering is as in \cite{miranda}, where the manifolds are called ``modes.'') 
Each computation is done using a single MCMC run of $n_t = 10^8$ points.
The error bars (such as $\pm0.12$ for manifold 8) are one standard deviation, 
estimated using \eqref{eq:sr}. 
We compare the calculations with the value $z_t$ obtained in \cite{miranda}.    }
\label{2dmodes}
\end{table}

\begin{figure}
\begin{center}
\begin{tikzpicture}

\node(1) at (0,0) {1};
\node(2) at (1,0) {2};
\node(3) at (2,0) {3};
\node(4) at (3,0) {4};
\node(5) at (4,0) {5};
\node(6) at (5,0) {6};

\draw (0,0) circle (0.5);
\draw (1,0) circle (0.5);
\draw (2,0) circle (0.5);
\draw (3,0) circle (0.5);
\draw (4,0) circle (0.5);
\draw (5,0) circle (0.5);

\node(1) at (8,0) {1};
\node(2) at (8.866,-0.5) {2};
\node(3) at (9.732,0) {3};
\node(4) at (9.732,1) {4};
\node(5) at (8.866,1.5) {5};
\node(6) at (8,1) {6};

\draw (8,0) circle (0.5);
\draw (8.866,-0.5) circle (0.5);
\draw (9.732,0) circle (0.5);
\draw (9.732,1) circle (0.5);
\draw (8.866,1.5) circle (0.5);
\draw (8,1) circle (0.5);

\end{tikzpicture}

\end{center}
\caption{Two dimensional view of two clusters of 6 spheres: a chain (left) and a loop (right).} 
\label{fig:worm_loop}
\end{figure}
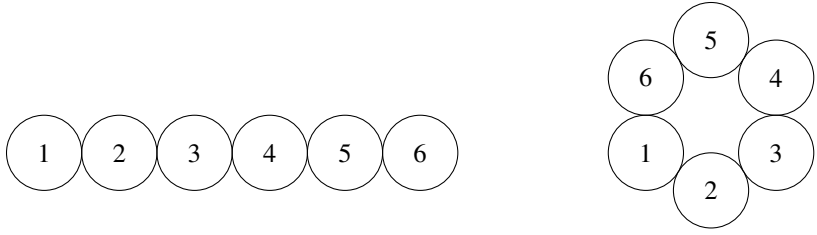

 Next, we consider a cluster of $N$ identical spheres which may form either a chain (C) or a loop (L), as in Figure \ref{fig:worm_loop}.  We are interested in the conditions under which chains are preferable to loops. 
 Let $M_{C,N}, M_{L,N}$ denote the manifolds corresponding to a chain, loop respectively; these have dimensions $2N-2$, $2N-3$. 
 For example, for $N=10$ spheres the dimensions are 18 and 17. 
 Integrals over these manifolds are not easily calculated by deterministic parameterizations. 
 
For each $N = 4-10$, and for each $i \in \{L,C\}$, we compute the volume $V$ and the geometrical partition function $z$ of the manifold $M_{i,N}$, using $n_t = 10^8$ points and $k = 4$ steps each. 
Because we don't have other results to compare to, we check that our calculations are correct by comparing our calculated value of the average of $f$, namely $\overline{h} = z/V$, to that estimated by our sampling algorithm, $\tilde h$. The values of $\tilde h$, $\bar h$ agree within statistical error as shown in Table \ref{LoopWorm}. 

With the values of $z_C, z_L$  we can make experimental predictions. Indeed, the ratio of the probabilities of finding a system of $N$ sticky sphere in a chain versus a loop  in equilibrium is 
 \begin{equation}\label{eq:LC}
 \frac{P(\text{chain})}{P(\text{loop})} = \kappa^{-1} \frac{n_C\,z_C}{n_L\,z_L}\,,
\end{equation}
where $n_C, n_L$ are the number of distinct copies of each manifold that one obtains by considering all permutations of spheres. If the spheres are indistinguishable then one can calculate\footnote{
The number of distinct copies is $N!/o$, where $o$ is the symmetry number of the clusters, i.e. the number of permutations that are equivalent to an overall rotation.
For the chain, the symmetry group is $C_2$, the cyclic group of order 2, and therefore the symmetry number is $2$. For the loop, the symmetry group is $C_N$, the rotational group of a regular N-gon of order $N$, so the symmetry number is $2N$. 
It is also possible that the spheres are distinguishable, for example if the bonds forming the chain are unbreakable and only one bond can break and reform, in which case we would have $n_C{=}n_L{=}1$. 

 Note that for  three-dimensional clusters, the number of distinct clusters of $N$ spheres should usually be computed as $2N!/o$, where $o$ is the number of permutations that are equivalent to a rotation of the cluster (e.g. \cite{wales}.) In our case, however, the clusters can be embedded into a two-dimensional plane, so a reflection of the cluster is equivalent to a rotation, hence we do not need to include the factor 2.
} 
that $n_C = N!/2$ and $n_L = (N{-}1)!/2$. 
From the calculated ratios $(n_Cz_C)/(n_Lz_L)$, shown in Table \ref{table:ratio}, and the value of the sticky parameter, which must be measured or inferred in some other way, we obtain a prediction for the fraction of chains versus loops, which could be compared to experimental data. Of course, the spheres could assemble into other configurations, whose entropies we could in principle calculate, but \eqref{eq:LC} still holds without calculating these additional entropies because the normalization factors cancel. 

 If we don't know the sticky parameter $\kappa$ for the system, then \eqref{eq:LC} can actually be used to \emph{measure} it for a given system. Indeed, an estimate for $\kappa$ using a fixed value of $N$ would be 
 \begin{equation}\label{eq:Kestimate}
\kappa \;\approx\; \widehat{\kappa} = \frac{\text{\# of loops of length $N$}}{\text{\# of chains of length N}}\cdot\frac{n_C\,z_C}{n_L\,z_L}\,.
 \end{equation}
Such a measurement could be useful, for example, in systems of emulsion droplets whose surfaces have been coated by strands of sticky DNA \cite{Feng:2013dr}. The DNA attaches the droplets together, like velcro, but the sticky patch can move around on the surface of a droplet, so a chain of attached droplets can sample many configurations. Understanding the strength of the droplet interaction from first principles, or even measuring it directly, is a challenge because it arises from the microscopic interactions of many weak DNA bonds as well as small elastic deformation of the droplet. But, \eqref{eq:LC} could be used to measure the stickiness simply by counting clusters in a microscope.

\begin{table}[!t]
\begin{center}
\begin{tabular}{l |  r @{\,} c@{\,} l   r @{\,} c@{\,} l  r@{.}l @{\,} c@{\,} r@{.}l   r@{.}l @{\,} c@{\,} l}
& \multicolumn{3}{c}{$V$}   &  \multicolumn{3}{c}{$z$} & \multicolumn{5}{c}{$\overline{h}$}   &  \multicolumn{4}{c}{$\widetilde{h}$}\\
\hline
4 Chain  & $3.51 \cdot 10^2$ & $\pm$ & 1.03 & $1.32 \cdot 10^2$ & $\pm$ & 0.38 & 0&374 & $\pm$ & 0&0022 & 0&379 & $\pm$ & $1.5\cdot 10^{-5}$  \\
4 Loop   & $8.27 \cdot 10^1$ & $\pm$ & 0.20 & $2.18 \cdot 10^1$ & $\pm$ & 0.054 & 0&264 & $\pm$ & 0&0013 & 0&261 & $\pm$ & $7.6\cdot 10^{-5}$  \\
5 Chain & $3.03 \cdot 10^3$ & $\pm$ & 16 & $8.31 \cdot 10^2$ & $\pm$ & 4.3 & 0&274 & $\pm$ & 0&0029 & 0&279 & $\pm$ & $1.1\cdot 10^{-4}$\\ 
5 Loop & $3.19 \cdot 10^2$ & $\pm$ & 4.3 & $6.19 \cdot 10^1$ & $\pm$ & 0.85 & 0&194 & $\pm$ & 0&0053 & 0&194 & $\pm$ & $2.8\cdot 10^{-4}$  \\
6 Chain & $2.70 \cdot 10^4$ & $\pm$ & $2.7 \cdot 10^2$  & $5.33 \cdot 10^3$ & $\pm$ & 52 & 0&197 & $\pm$ & 0&0039  & 0&206 & $\pm$ & $1.9\cdot 10^{-4}$  \\
6 Loop & $1.66 \cdot 10^3$ & $\pm$ & 13 & $2.32 \cdot 10^2$ & $\pm$ & 1.9 & 0&143 & $\pm$ & 0&0024 & 0&147 & $\pm$ & $3.1\cdot 10^{-5}$  \\
7 Chain & $2.42 \cdot 10^5$ & $\pm$ & $1.3\cdot 10^4$ & $3.84\cdot 10^4$ & $\pm$ & $1.9\cdot 10^3$ & 0&158 & $\pm$ & 0&0017 & 0&151 & $\pm$ & $1.7\cdot 10^{-4}$\\ 
7 Loop & $1.01 \cdot 10^4$ & $\pm$ & $1.2\cdot 10^2$ & $1.09\cdot 10^3$ & $\pm$ & $13$ & 0&109 & $\pm$ & 0&0026 & 0&110 & $\pm$ & $2.6\cdot 10^{-5}$\\ 
8 Chain & $2.19 \cdot 10^6$ & $\pm$ & $6.4\cdot 10^4$ & $2.77\cdot 10^5$ & $\pm$ & $8.3 \cdot 10^3$ & 0&127 & $\pm$ & 0&0075 & 0&111 & $\pm$ & $9.6\cdot 10^{-4}$\\ 
8 Loop & $6.76 \cdot 10^4$ & $\pm$ & $1.6\cdot 10^3$ & $5.65\cdot 10^3$ & $\pm$ & $1.3 \cdot 10^2$ & 0&0836 & $\pm$ & 0&0039 & 0&0824 & $\pm$ & $2.5\cdot 10^{-5}$\\ 
9 Chain  & $2.12 \cdot 10^7$ & $\pm$ & $1.1\cdot 10^6$ & $1.73\cdot 10^6$ & $\pm$ & $8.9 \cdot 10^4$ & 0&0816 & $\pm$ & 0&0084 & 0&0820 & $\pm$ & $1.1\cdot 10^{-5}$\\ 
9 Loop & $4.94 \cdot 10^5$ & $\pm$ & $1.7\cdot 10^4$ & $2.87\cdot 10^4$ & $\pm$ & $9.2 \cdot 10^2$ & 0&0582 & $\pm$ & 0&0038 & 0&0612 & $\pm$ & $1.2\cdot 10^{-5}$\\ 
10 Chain & $2.07 \cdot 10^8$ & $\pm$ & $1.2\cdot 10^7$ & $1.43\cdot 10^7$ & $\pm$ & $8.1 \cdot 10^5$ & 0&0689 & $\pm$ & 0&0080 & 0&0604 & $\pm$ & $3.5\cdot 10^{-4}$\\ 
10 Loop & $3.81 \cdot 10^6$ & $\pm$ & $2.2\cdot 10^5$ & $1.71\cdot 10^5$ & $\pm$ & $9.7 \cdot 10^3$ & 0&0448 & $\pm$ & 0&0051 & 0&0453 & $\pm$ & $4.7\cdot 10^{-5}$
\end{tabular}
\end{center}
\caption{Volume $V$, geometrical partition function $z$ and mean vibrational contribution $\overline{h} = z/V$ for the manifolds $M_{i,N}$ corresponding to the loop/chains of $N$ particles. The computations were done using $n_t = 10^8$ points. Error bars are calculated using \eqref{eq:sr}. The results are compared with the average vibrational contribution $\widetilde{h}$ estimated from our sampling algorithm. 
Specifically, we generate $n_p$ points on  $M_{i,N}$ (and repeat this $n_m=10$ times) and estimate the average as $\widetilde{h} = \frac{1}{n_m}\sum_{j=1}^{n_m} \widetilde{h}_j$ where $\widetilde{h}_j = \frac{1}{n_p} \sum_{i=1}^{n_p} f(X_i)$. 
The error bars of the latter are computed using the sample standard deviation, using a sample of size $n_m = 10$.}
\label{LoopWorm}
\end{table}
 
\begin{table}[!t]
\begin{center}
\begin{tabular}{c r@{.}l @{\,} c@{\,} l   r@{}l @{\,} c@{\,} r@{}l  }
$N$ &  \multicolumn{4}{c}{$z_C/z_L$}   &  \multicolumn{5}{c}{$(n_Lz_L)/(n_Cz_C)$}\\
\hline
$4$ & 6&02 &$\pm$& 0.032  & 24.&1 &$\pm$& 0.&13 \\
$5$ & 13&5 &$\pm$& 0.13  & 67.&5 &$\pm$& 0.&65  \\
$6$ & 22&9 &$\pm$& 0.41 & 138& &$\pm$& 2.&5  \\
$7$ & 35&1 &$\pm$& 1.0  & 246& &$\pm$& 7.&0 \\
$8$ & 49&0 &$\pm$& 2.6  & 392& &$\pm$& 21& \\
$9$ & 60&3 &$\pm$& 5.0  & 543& &$\pm$& 45& \\
$10$ & 83&6 &$\pm$& 9.5  & 836& &$\pm$& 95& \\
\end{tabular}
\end{center}
\caption{Ratios between the geometrical partition functions of a chain and a loop of $N$ particles, for $4\leq N\leq10$. We report both the ratio $z_C/z_L$ of a single chain vs. loop, and the ratio $(n_Cz_C)/(n_Lz_L)$ when spheres are treated as indistinguishable. The latter may be used to estimate the sticky parameter $\kappa$ as in \eqref{eq:Kestimate}. 
}
\label{table:ratio}
\end{table}

 Finally, we note that the ratio $(n_C z_C)/(n_Lz_L)$ is related to the loss of entropy as a chain forms a bond to become a loop. Even though a loop has more bonds, hence lower energy, it still might be less likely to be observed in equilibrium because it has lower entropy. The ratio $(n_C z_C)/(n_Lz_L)$ is exactly the value of stickiness $\kappa$ above which a loop becomes more likely than a chain, which corresponds to a threshold temperature or binding energy. 
 Note that as $N$ increases the ratio does too, so we expect fewer loops for longer polymers than for shorter ones.

\section{Minimizing the error bars: some toy models to understand parameter dependence}\label{sec:toy}

Equation \eqref{eq:sr} estimates the relative standard deviation of the estimator $\widehat{Z}$ for $Z$. 
We would like to choose parameters in the algorithm to minimize this relative error. Previous groups studying similar algorithms in Euclidean spaces have investigated how the overall number of operations scales with the dimension \cite[e.g.][]{vempala}, but to our knowledge there is little analysis of how such an algorithm behaves for fixed dimension, and how to optimize its performance with fixed computing power. 

To this aim we construct a hierarchy of toy models that suggest how the error formula may depend 
on parameters. We focus exclusively on the sum $\sum_i (1-p_i)\tau_i / (n_ip_i)$ at the right of \eqref{eq:sr} which estimates the relative variance of a product of ratio estimates. This is because the estimator $\widehat{Z}_k$ for $\rho_k$ usually has a much smaller error, since it is estimated by straightforward Monte Carlo integration with independently chosen points. 
Our models may apply to other multi-phase MCMC algorithms.  

We consider an algorithm which chooses the same number of points $n$ at each stage and has a fixed number of total points $n_t = kn$, where $k$ is the number of stages. 
Additionally, we fix the largest and smallest radii $r_0$, $r_k$, and take the $d$-dimensional volume ratios between balls to be the same constant $\nu_i = \nu$ at each stage. The ratio $\nu$ and number of stages $k$ are related by 
\begin{equation}\label{eq:nuk}
\nu^k = \left(\frac{r_0}{r_k}\right)^d = C, 
\end{equation}
where we will write $C$ to denote a constant (not always the same) which does not depend on $\nu$ or $k$.

We now make the approximation that $R_i \propto\nu$. This would be exactly true if $M$ were a hyperplane of dimension $d$. 
The relative variance we wish to minimize is
\begin{equation}  \label{eq:g}
     g(\nu) = C \frac{\nu-1}{\log(\nu)} \sum_{i=0}^{k-1} \tau_i \; .
\end{equation}
 Recall that $\tau_i$ is the correlation time for the indicator function  ${\bf 1}_{B_{i+1}}(X_j)$, assuming points are generated in $M\cap B_i$ (see \eqref{eq:tau}, \eqref{eq:vNi}.) 
 
 We now consider several models for the correlation times $\tau_i$ that will allow us to investigate the behavior of $g(\nu)$. Ultimately, we would like to find the value $\nu^*$ which minimizes this function.

\subsection{Constant $\tau$.}

Suppose the $\tau_i= \tau$ are equal and independent of $k$.
This would happen, for example, if we were able to make independent samples and achieve $\tau_i = 1$.
It also could happen if there were a scale invariant MCMC algorithm such as the hit-and-run
method or affine invariant samplers.
It seems unlikely that such samplers exists for sampling embedded manifolds.
In that case, we have
\begin{equation}
       g(\nu) = C \frac{\nu-1}{\log(\nu)^2} \; .
\end{equation}
Optimizing this over $\nu$ gives $\nu^* \approx 4.9$, independent of dimension.
This result is interesting in that it suggests a more aggressive strategy than
the previously-suggested choice $\nu = 2$ \cite{vempala}, i.e. one should shrink the balls by a larger amount and hence have a smaller number of stages.

\subsection{Diffusive scaling}\label{sec:diffusive}

As a second scenario, we model $\tau_i$ as 
\begin{equation}\label{eq:tau_model}
   \tau_i \propto r_i^2.  
\end{equation}
This is a ``diffusive scaling'' because it represents the time scale of a simple diffusion or random
walk in a ball of radius $r_i$ (see \eqref{hpn} in section \ref{sec:brownian_motion} below.)
From \eqref{eq:nu}, we have $\tau_i \propto \nu^{-2i/d}$. 
Hence \eqref{eq:g} becomes, after some algebra and using \eqref{eq:nuk}:
\begin{equation}\label{eq:var_approx_fin}
  g(\nu) \propto g_d(\nu) \equiv \frac{\nu-1}{\log(\nu)(1-\nu^{-2/d})} \;.
\end{equation}

\begin{figure}[!t]
\centering
\includegraphics[scale = 0.25]{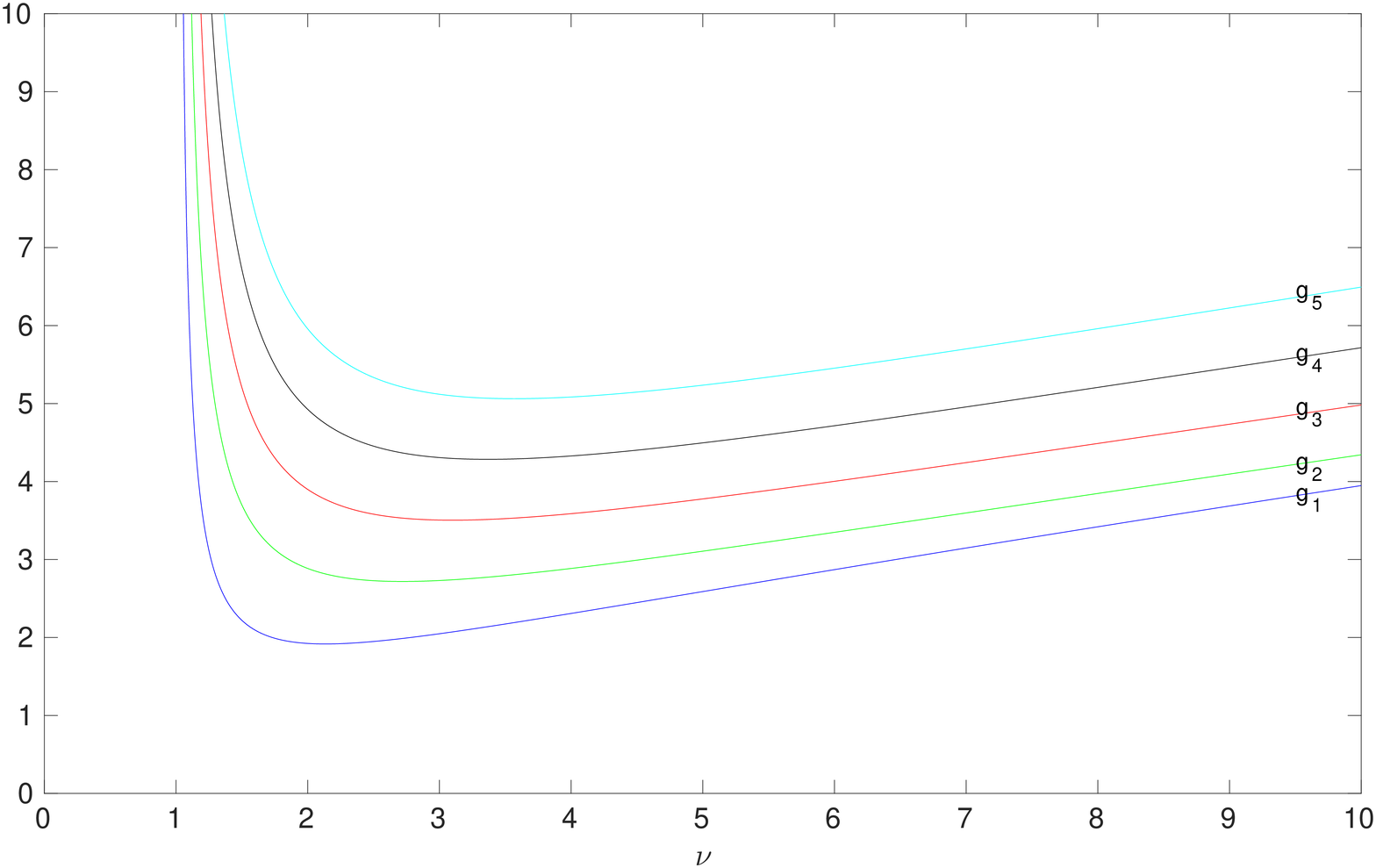}
\caption{Plots of the functions $g_d(\nu)$ as in \eqref{eq:var_approx_fin} for $d = 1,2,3,4,5$.
}
\label{fig:variance_r2}
\end{figure}

Plots of $g_d(\nu)$ are given in Figure \ref{fig:variance_r2}.
Examples of minimizers $\nu_d^*$ are 
$\{2.6,2.7,3.1,3.4,3.6,3.7,4.1,4.5,4.7\}$ for $d=\{1,2,3,4,5,6,10,20,50\}$ respectively. 
Some analysis shows that $\nu^*_d \to \nu^*$ as $d \to \infty$, where $\nu^*$ is the optimal
ratio assuming constant $\tau$ above.
The same analysis shows that $g_d(\nu^*_d)\approx \frac{d(\nu^*-1)}{2\log^2(\nu^*)}$ 
for large $d$, so the relative variance in this model is proportional to the manifold's intrinsic dimension $d$. 
The figure shows that $g_d(\nu)$ increases more slowly for $\nu > \nu^*$ than for $\nu < \nu^*$.
This also suggests taking larger values of $\nu$ in computations.

We found it surprising that the model with diffusive scaling predicted smaller ratios than the one with constant scaling, since it has a correlation time which decreases as the balls shrink so one might prefer to shrink them rapidly. Some insight into why comes from thinking about the \emph{marginal} cost of adding an additional ball, given we already have $k'$ balls. For the first model, this marginal cost is constant no matter what $k'$ is, and equal to the cost of just the first ball. For the second model, the marginal cost can decrease as $k'$ increases, since the contribution to the cost from the $j$th ball is a decreasing function of  $j$. So, with diffusive scaling, each additional ball does not increase the error by as much. 
Another way of seeing this is to plot the correlation times as a function of ball number, and note that the first and last balls have fixed correlation times. The total error is roughly the area under this plot. 
But, the area under an exponential function connecting $(0,Cr_0^d)$ to $(k',Cr_k^d)$ increases less, when $k'\to k'+1$, than the area under a constant function connecting $(0,C')$ to $(k',C')$ (here $C$, $C'$ are constants.)

\subsection{A Brownian motion model in balls.}\label{sec:brownian_motion}

As a third model, we approximate the MCMC 
sequence $X_0,X_1,\ldots$ with a continuous-time Brownian motion, and simplify the geometry to two
concentric balls in a $d-$dimensional plane. 
This allows us to find an analytic formula for $\tau_i$. 
While this is a natural model to consider, we will describe later why it fails to capture the behavior of almost any system of interest. Nevertheless, the reasons behind its failure help build intuition into high-dimensional geometry, so we feel it useful to describe this mathematical exercise. 

What follows is an exact derivation of the correlation times in this model. Readers interested in only the result may skip to equation \eqref{eq:var_balls}  and the subsequent discussion. 

Let $X_t$ be a $d$-dimensional Brownian motion on a $d$-dimensional ball $B_0$ of radius $R_0$
centered at the origin, with reflecting boundary conditions at $\left|x\right|=R_0$.
We consider the stochastic process 
\begin{equation}\label{eq:indic}
    F(X_t) = {\bf 1}_{B_1}(X_t), 
\end{equation}
where $B_1$ is a ball with radius $R_1 < R_0$, also centered at the origin. 
Our goal is to derive an analytic formula for the correlation time of the process, defined for a continuous-time process to be (compare with (\ref{eq:tau})) \cite{sokal}
\begin{equation}\label{tau}
   \tau = 2 \int_0^{+\infty} \frac{C(t)}{C(0)}\, dt \; .
\end{equation} 
This is expressed in terms of the 
stationary covariance function of $F(X_t)$: 
\begin{equation}\label{eq:C}
C(t) = \mathbb{E}\!\left[\,({\bf 1}_{B_1}(X_0) - \overline{F})
                              ({\bf 1}_{B_1}(X_t) - \overline{F})\right], \qquad
 \overline{F} = \mathbb{E}\!\left[ {\bf 1}_{B_1}(X)\right] .
\end{equation}
The expectation is calculated assuming that $X_0$ is drawn from the stationary distribution, which for a Brownian motion is the uniform distribution. 

We know that 
\begin{equation}\label{eq:C0}
C(0) = \mbox{var}(F(X_t)) = \overline{F}(1-\overline{F})
\end{equation}
 (as for a Bernoulli random variable), so we focus on computing the integral in \eqref{tau}.  
 Notice that $\overline{F} = (R_1/R_0)^d$. 

We start by writing \eqref{eq:C} in terms of the probability density $p(t,x,y)$ for $X_t=x$ given $X_0=y$.
If $X_0$ is in the stationary distribution, then the joint density of $X_0=y$ and $X_t = x$ 
is $Z_0^{-1} p(t,x,y)$, where $Z_{0} = \mbox{vol}(B_0)$.  
Using this expression to write
\eqref{eq:C} as an integral over $dx,dy$, integrating over time, and rearranging the integrals gives 
\begin{equation} \label{eq:C(t)}
  \int_0^\infty   C(t) dt = Z_0^{-1}\int_{B_0} \int_{B_0}({\bf 1}_{B_1}(x) - \overline{F})
                                  ({\bf 1}_{B_1}(y) - \overline{F})\int_0^\infty p(t,x,y) \,dtdydx \; .
\end{equation}
To calculate the time integral we must evaluate
$
\overline{p}(x,y) = \int_0^{+\infty} p(t,x,y) \,dt .
$
Notice that $p(t,x,y)$ satisfies the forward Kolmogorov  
equation with boundary and initial conditions (see, e.g., \cite{gardiner}):
\begin{equation}\label{fp}
p_t = \frac{1}{2} \Delta_x p \;\;  (x \in B_0) , \quad
\nabla p \cdot \hat{\bfn} = 0 \;\;  (x \in \partial B_0) , \quad
 p(0,x,y) = \delta(x-y) ,
\end{equation}  
where $\hat{\bfn}$ is the vector normal to the boundary. 
Therefore $\overline{p}(x,y)$ satisfies the equation
\begin{align*}
    \Delta_x \overline{p}(x,y) = \int_0^{+\infty} \Delta_x p(t,x,y) \,dt  
                               =  2 \int_0^{+\infty} p_t(t,x,y) \mathrm{d}t 
                               = -2 \delta(x-y)\; .
\end{align*}
Therefore $\overline{p}(x,y)$ is proportional to the Green's function for the Laplace equation in a ball.
We now wish to evaluate the integral over the $y$ coordinate
\[
u(x)  = \int_{B_0} f(y)\overline{p}(x,y) \mathrm{d}y ,
\]
where we have defined
$f(x) = {\bf 1}_{B_1}(x) - \overline{F}$ to be the centered function whose covariance we wish to evaluate. 
One can check that $u(x)$ solves the Poisson equation 
\begin{equation}\label{u}
 \Delta_x u = -2 f \quad \text{for } x \in B_0, \qquad
 \nabla u \cdot \hat{\bfn} = 0 \quad \text{for }x \in \partial B_0 . 
\end{equation} 
We now solve this equation analytically for $u(x)$. 
For $d =1$, \eqref{u} reduces to an ODE, whose solution is 
\begin{equation}\label{sol1}
u(x) = \left\{ \begin{aligned}
\left( \frac{R_1}{R_0}-1 \right) x^2 & \quad |x| < R_1 \\
\frac{R_1}{R_0} x^2-2R_1|x|+R_1^2 & \quad R_1 < |x| < R_0
\end{aligned} \right.
\end{equation}

For $d \geq 2$, we solve \eqref{u} by looking at radial solutions $u = u(|x|) = u(r)$. 
In polar coordinates \eqref{u} is
\begin{equation}\label{upolar}
 u_{rr} + \frac{d-1}{r}u_r  = -2 f \quad\text{for } r \in (0,R_0), \qquad 
 u_r(R_0) = 0.
\end{equation} 
We solve this equation using the elementary radial solutions  $\Phi(r)$ of the Laplace equation ($\Delta u = 0$) in the disk,\footnote{These are given by 
\begin{equation*}
\Phi(r) = \left\{ \begin{aligned}
& C \; \text{log}(r) + K & d = 2 \\
& \frac{C}{2-d}r^{2-d}+ K & d \geq 3
\end{aligned} \right.
\end{equation*} 
where $C$ and $K$ are constants \cite{evans}.
}
and requiring the solution to satisfy $\lim_{r \to 0} u(r) = 0$ and to be continuous at $r = R_1$. 
For $d = 2$ we have
\begin{equation}\label{usol2}
u(r) = \left\{ \begin{aligned}
& \frac{R_1^2-R_0^2}{R_0^2} r^2 & \quad 0\leq r\leq R_1 \\
& R_1^2(\text{log}R_1-\text{log}r)+\frac{R_1^2}{2}\left(\frac{r^2}{R_0^2}-1\right) & R_1 < r \leq R_0
\end{aligned}  \right.
\end{equation}
and for $d \geq 3$ we have
\begin{equation}\label{usol3}
u(r) = \left\{
\begin{aligned}
& \frac{R_1^d-R_0^d}{dR_0^d} r^2 & 0 \leq r \leq R_1 \\
& -\frac{2R_1^d}{d(2-d)} r^{2-d} + \frac{R_1^d}{dR_0^d}r^2 + \frac{R_1^2}{2-d} & R_1 \leq r \leq R_0
\end{aligned} \right.
\end{equation}

We can now evaluate the integral over the $x$ coordinate to obtain 
\begin{equation}\label{eq:idx}
   \int_0^{+\infty} C(t) \,dt =Z_0^{-1} \int_{B_0} f(x) u(x) \, dx = \mbox{vol}(S_{d-1})Z_0^{-1}\int_0^{R_0} f(r)u(r)r^{d-1} \mathrm{d}r 
\end{equation}
where in the last step we have written the integral using polar coordinates. 
Here $\mbox{vol}(S_{d-1})$ is the surface area of the boundary of the unit $d$-dimensional sphere. 
Substituting our analytic expression for $u(x)$, and using \eqref{eq:C0} in \eqref{tau}, we obtain
\begin{equation}\label{hpn}
\tau = R_0^2 h_d(\nu) \;.
\end{equation}
Here $\nu = (R_o/R_1)^d$ is the volume ratio as in \eqref{eq:nu}, and the functions $h_d(\nu)$ are given by
\begin{equation}\label{tau_sol}
h_d(\nu) = \left\{ \begin{aligned}
& \frac{4}{3} \frac{\nu-1}{\nu^2}  & d = 1 \\
&  \frac{\nu^{-1}-1+\text{log}(\nu)}{\nu-1} & d = 2 \\
& \frac{4}{d^2-4} \frac{(d-2)\nu^{-1} -d\nu^{2/d-1}+2}{\nu^{2/d}(1-\nu^{-1})}& d \geq 3
\end{aligned}
\right.
\end{equation}

\begin{figure}[!t]
\begin{center}
\includegraphics[scale = 0.28]{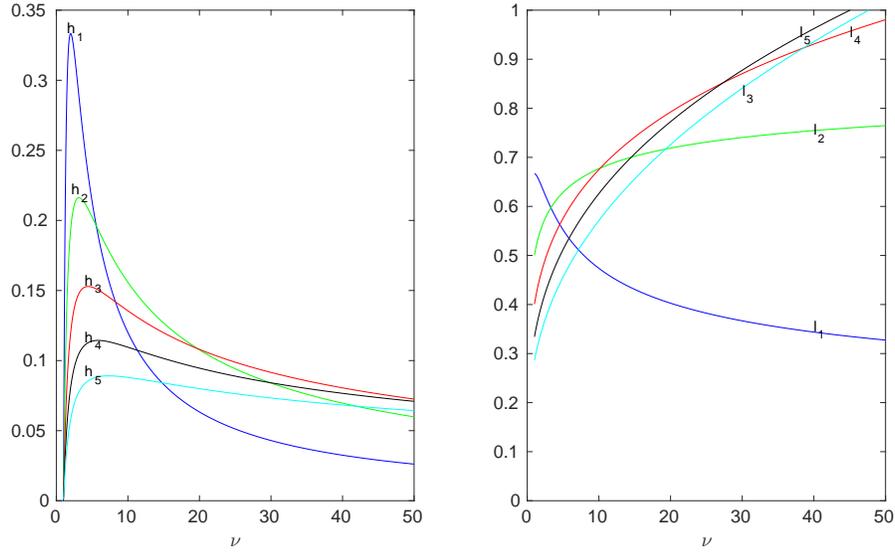}
\caption{Left: plots of $h_d(\nu)$ from \eqref{tau_sol}, which are proportional to the correlation times of the process ${\bf 1}_{B_1}(X_t) $ in \eqref{eq:indic}.
Right: plots of $l_d(\nu)$ in  \eqref{eq:var_balls}, which are proportional to the relative variance of the product of ratios. Both functions are analytically computed for a Brownian motion $X_t$ in two concentric spheres, for dimensions $d = 1,2,3,4,5$.
}
\label{fig:variance_balls}
\end{center}
\end{figure}

The functions $h_d(\nu)$ are proportional to the correlation time and are plotted in Figure \ref{fig:variance_balls}  for $d = 1,2,3,4,5$. For fixed $\nu$ and $d\geq 3$ the correlation time decreases with the dimension, and for all $d$ we have 
\begin{equation*}
\lim_{\nu \to 1} h_d(\nu) = \lim_{\nu \to +\infty} h_d(\nu) = 0.
\end{equation*}
We found all of these observations to be surprising and will discuss them further in section \ref{sec:compare}.

We use \eqref{hpn} to approximate the correlation time $\tau_i$ in \eqref{eq:g} as 
\begin{equation}\label{eq:tau_hdn}
\tau_i = r_{i}^2 h_d(\nu).
\end{equation} 
After some algebra, we obtain
\begin{equation}\label{eq:var_balls} 
g(\nu) \propto l_d(\nu) \equiv \frac{(\nu-1)h_d(\nu)}{\text{log}(\nu)(1-\nu^{-2/d})}
\end{equation}
One can verify from \eqref{eq:var_approx_fin}, \eqref{eq:var_balls} that $l_d(\nu) = g_d(\nu)h_d(\nu)$, so this model's only difference with that of section \ref{sec:diffusive} is that it additionally models how $\tau_i$ depends on $\nu$. 

Plots of the function $l_d(\nu)$ for dimension up to $d = 5$ are given in Figure \ref{fig:variance_balls}. The function $l_1(\nu)$ is strictly decreasing, while $l_d(\nu)$ is strictly increasing for $d\geq2$. It follows that the ratio $\nu^*$ that minimizes \eqref{eq:var_balls} is $\nu^* = 1$ for $d \geq 2$. 

\subsection{Comparing the toy models with data}\label{sec:compare}

\begin{figure}[!t]
\center
\includegraphics[scale = 0.3]{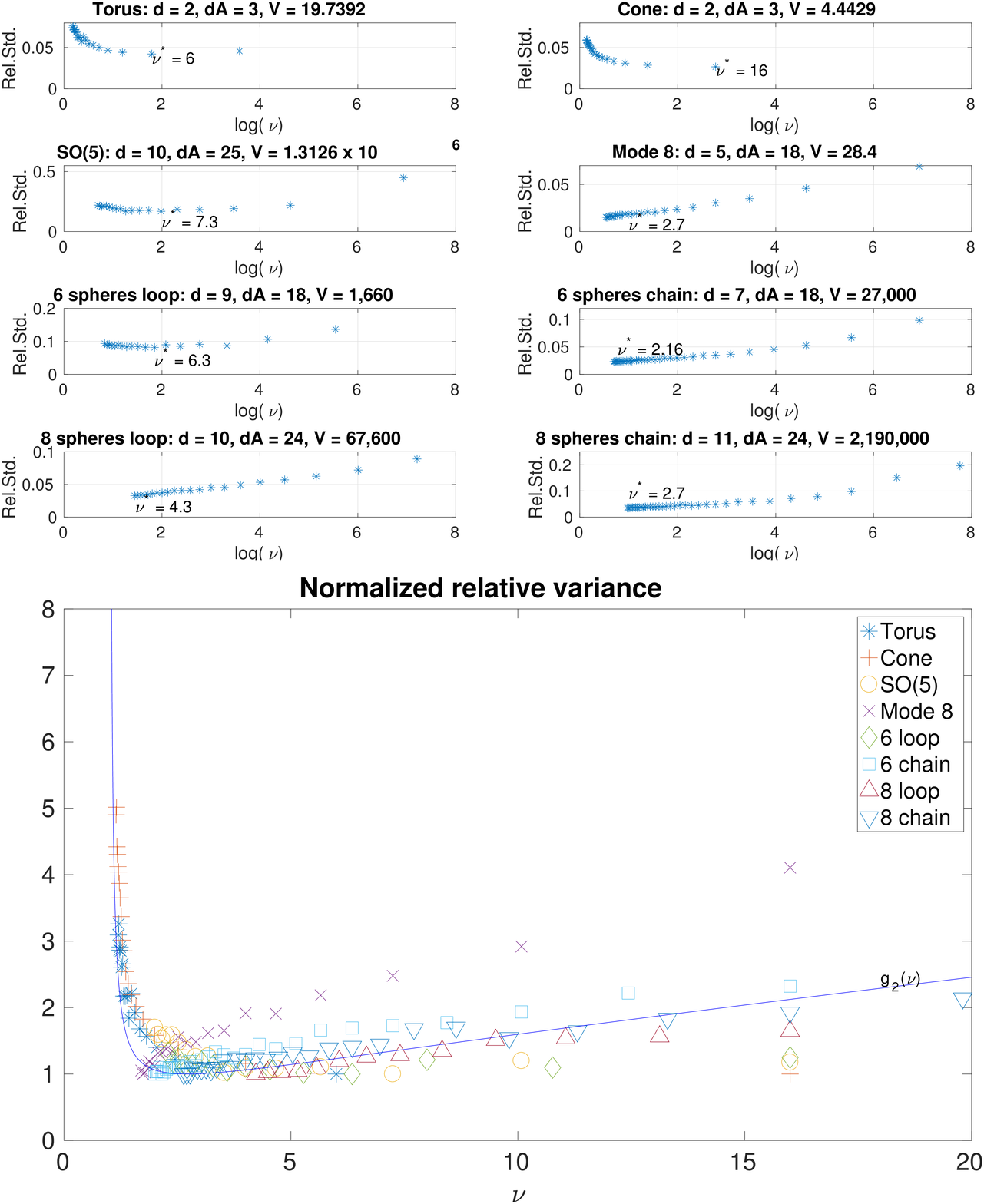}
\caption{Top: plots of the relative standard deviation $\sigma_r$ of $\widehat{Z}$ as a function of the ball ratio $\nu$ for eight different manifolds. These are each obtained using formula \eqref{eq:sr} from a single run of $n_t = 10^6$ points. The minimizers $\nu^*$ are reported in the plots and plotted as smaller black dots. Bottom: plots of the relative variance $\sigma^2_r$ of $\widehat{Z}$ of the same manifolds, combined together and normalized so that they all have the minimum in $\nu = 1$, superimposed on the theoretical model $g_2(\nu)$ given in \eqref{eq:var_approx_fin}. }
\label{fig:errors}
\end{figure}

Our first two models for the correlation time 
predicted an error that was minimized at an intermediate value $\nu^*$, equal to about $\nu^*\approx 5$ for the model with constant $\tau$, and slightly smaller for the model with diffusive scaling though approaching it as $d\to \infty$.
The third was totally different -- it predicted taking as small  $\nu$ as possible. What does the empirical evidence suggest? 

Figure \ref{fig:errors} plots the relative standard deviation of $\widehat{Z}$ obtained empirically from averaging over several simulations, as a function of $\log(\nu)$, for several different manifolds of varying dimensions and co-dimensions. 
For this figure we calculated the volume using $k = 1, \ldots, 40$ and computed $\nu$ using \eqref{eq:nuk}, assuming fixed $r_0,r_k$. 
All curves have the same basic shape, which is qualitatively the same as that predicted by the first two models but not the third.

What went wrong with the third model? The only difference it has with the diffusive scaling model is an additional assumption for how the correlation time depends on $\nu$. 
This correlation time was small for large $\nu$, something to be expected since a Brownian motion  should leave $B_1$ very quickly when $B_1$ is much smaller than $B_0$. 
But, importantly, the correlation time was also small for small $\nu$; small enough to make the overall relative variance decrease monotonically as $\nu$ decreases. 
The fact that the correlation time approaches 0 as $\nu\to 1$ 
implies that a Brownian motion leaves $B_1$ very quickly even when $B_1$ is nearly as big as $B_0$.

We believe this is an artifact of our highly restrictive geometry, and almost any problem of interest will not have this property. In high dimensions, virtually all of the volume of an object will be concentrated near the extremal parts of its boundary, a fact that can be checked for example by calculating the volume of a spherical shell of width $\epsilon\ll 1$. 
A sphere's boundary therefore is concentrated near the outermost spherical shell, so when a Brownian motion enters $B_1$ it will with overwhelming probability stay near the boundary, where it is easy for it to leave $B_1$. Therefore, we expect the Brownian motion to quickly enter and leave $B_1$ repeatedly, like a Brownian motion jumping back and forth between the boundaries of a thin strip. The correlation time therefore should be short.

When $B_1\cap M$, $B_0\cap M$ are not spherical, then some parts of their boundaries will lie further from the origin than others, and it is in these ``corners'' that we expect a Brownian motion to get stuck for longer times. In convex geometry such corners are known to control the computational 
complexity of Monte-Carlo algorithms \cite{simonovits}, an effect that is simply not possible to capture with Brownian motion in spheres. We expect that inhomogeneous boundary densities are the norm rather than the exception in any problem of scientific interest, and therefore do not believe the predictions of this third model.


\section{Conclusions}  \label{sec:conclusions} 


We introduced algorithms to sample and integrate over a  connected component of a manifold $M$, defined by equality and 
inequality constraints in Euclidean space. 
Our sampler is a random walk on $M$ that does not require an explicit parameterization of 
$M$ and uses only the first derivatives of the equality constraint functions. 
This is simpler than other samplers for manifolds that use second derivatives.
Our integration algorithm was adapted from the multi-phase Monte Carlo methods for computing high dimensional volumes, based on intersecting the manifold with a decreasing sequence of balls.      

We discussed several toy models aimed at understanding how the integration errors depend on the factor by which the balls are shrunk at each step. Our models and data suggested taking larger factors than those suggested in the literature. 
We also estimated these errors by computing results analytically for a Brownian motion in balls,
but the results of these calculations disagree with computational experiments with the sampler.
We discussed the geometric intuition behind this failure. 

We tested the algorithms on several manifolds of different dimensions whose distributions or volumes are known analytically, and then we used our methods to compute the entropies of various clusters of hard sticky spheres. These calculations make specific predictions that could be tested experimentally, for example using colloids interacting via sticky DNA strands. 

We expect our methods to apply to other systems of objects that can be effectively modeled as ``sticky,'' i.e. interacting with a short-range potential \cite{colloids3}. For example, one could calculate the entropies of self-assembling polyhedra, a system that has been realized experimentally in the hope that it will lead to efficient drug-delivery methods \cite{Pandey:2011jj}; one could calculate the entire free energy landscape of small clusters (building on results in \cite{miranda}); or one could more accurately calculate the phase diagrams of colloidal crystals \cite{Macfarlane:2011fh,Tkachenko:2016ch} by computing entropic contributions to the free energy.

While our methods have worked well on the systems considered in this paper, we expect that applying them to certain larger systems will require additional techniques and ideas. For example, the manifold corresponding to the floppy degrees of freedom of a body-centered cubic crystal is thought to have a shape like a bicycle wheel, with long thin spokes leading out of a small round hub near the center of the crystal \cite{Jenkins:2014js}. Sampling efficiently from a manifold with such different length scales will be a challenge, though one may draw on techniques developed to sample the basins of attraction of jammed sphere packings, for which the same issues arise \cite{Martiniani:2016bt}.   
A different, though possibly related problem, is to sample from manifolds whose constraint functions are not linearly independent. Such manifolds arise for example in clusters of identical spheres \cite{Kallus:2016uq} or origami-based models \cite{Johnson2016}, and could be important for understanding how crystals nucleate.

There are many possible improvements to the Monte Carlo methods described here.
For the sampling algorithm, one improvement would be to use gradient information of the density function $f$ to create a more effective proposal distribution.  
It may be possible to find a manifold version of exact discrete Langevin methods
described, for example, in \cite{sokal} and \cite{liu}.
An interesting extension to the sampler would be to allow it to jump between manifolds of different dimensions. With such a sampler, one could calculate the ratio of entropies in \eqref{eq:LC} directly, without first calculating the volumes. 

Our integration algorithm could also be enhanced in several ways. 
It could be extended to the case when the intermediate manifolds $M\cap B_i$ are not connected, for example by using umbrella sampling techniques.
It could be more adaptive: the radii $R_i$ could be chosen more systematically using variance estimators, as is
done in thermodynamic integration \cite{hou}. 
We could use additional ratios such as $N_{i,i+2}$ (see \eqref{nij}) to improve the ratio estimates. 
We could also let $n_i$ vary with ball radius so the errors at each ratio more equally distributed; according to the diffusive scaling model of section \ref{sec:diffusive}, we should use more points to estimate the larger ratios. 
Finally, it might be possible to change the probability distribution using smoother ``steps'' when breaking up the volume into a product of ratios, as in \cite{hou} and thermodynamic integration \cite{frenkel84}, though how to adapt this strategy to manifolds is not immediately clear.

\subsection*{Acknowledgements}
E. Z. and M. H.-C. acknowledge support from DOE grant DE-SC0012296. 
The authors would additionally like to thank Gabriel Stoltz and Tony Lelievre for interesting discussions that helped to improve this article.

\bibliographystyle{unsrt}
\bibliography{mcmc_biblio}

\end{document}